\documentclass[leqno,11pt]{article}
\usepackage{amssymb,latexsym,mathrsfs,color}
\usepackage[dvips]{graphicx} 
\usepackage{url}

 \catcode`@=11\def\@oddhead{\hbox{}\hfil\rm\thepage}\def\@oddfoot{}
 \def\@evenhead{\hbox{}\hfil\rm\thepage}\def\@evenfoot{}
 \@twosidetrue \catcode`@=12
\newtheorem{prp}{Proposition}
\newtheorem{lem}[prp]{Lemma}\newtheorem{thm}[prp]{Theorem}
\newtheorem{cor}[prp]{Corollary}
\newenvironment{prf}{\begin{trivlist}\item[\emph{Proof.}]}{\end{trivlist}
  \medskip\par}
\newenvironment{prfof}[1]{\begin{trivlist}\item[\emph{Proof of #1.}]}{
  \end{trivlist} \medskip \par}
\newenvironment{rem}{\begin{trivlist}\item[\emph{Remark.}]}{\end{trivlist}
  \medskip\par}
\def\prpb{\begin{prp}}\def\prpe{\end{prp}}
\def\lemb{\begin{lem}}\def\leme{\end{lem}}
\def\thmb{\begin{thm}}\def\thme{\end{thm}}
\def\corb{\begin{cor}}\def\core{\end{cor}}
\def\prfb{\begin{prf}}\def\prfe{\end{prf}}
\def\prfofb#1{\begin{prfof}{#1}}\def\prfofe{\end{prfof}}
\def\remb{\begin{rem}}\def\reme{\end{rem}}
\def\prpa#1{\label{p:#1}}
\def\lema#1{\label{l:#1}}\def\lemu#1{Lemma~\ref{l:#1}}
\def\thma#1{\label{t:#1}}\def\thmu#1{Theorem~\ref{t:#1}}
\def\cora#1{\label{c:#1}}\def\coru#1{Corollary~\ref{c:#1}}
\def\seca#1{\label{s:#1}}\def\secu#1{Section~\ref{s:#1}}

\def\itmb{\begin{enumerate}}\def\itme{\end{enumerate}}
\def\itdb{\begin{itemize}}\def\itde{\end{itemize}}
\def\ittb{\begin{description}}\def\itte{\end{description}}
\def\eqnb{\begin{equation}}\def\eqne{\end{equation}}
\def\arrb#1{\begin{array}{#1}}\def\arre{\end{array}}
\def\tabb#1{\par\noindent\begin{tabular}{#1}}
\def\tabe{\end{tabular}\par\noindent}
\def\eqna#1{\label{e:#1}}\def\eqnu#1{(\ref{e:#1})}

\def\QED{\relax\ifmmode\let\@tempa\relax\ifcase\@eqcnt\def\@tempa{& & &}\or
  \def\@tempa{& &}\else\def\@tempa{&}\fi\@tempa $\Box$ \else\hfill $\Box$ \fi}
\def\DDD{\relax\ifmmode\let\@tempa\relax\ifcase\@eqcnt\def\@tempa{& & &}\or
 \def\@tempa{& &}\else\def\@tempa{&}\fi\@tempa $\Diamond$
 \else\hfill $\Diamond$ \fi}

\def\Rom#1{\uppercase\expandafter{\romannumeral#1}}
\def\dsp{\displaystyle}

\def\limf#1{\displaystyle \lim_{#1\to\infty}}

\def\norm#1{\dsp\left\| #1 \right\|}

\def\Ccomb#1#2{\setbox0=\hbox{$\displaystyle\mathrm{C}$}\setbox1=\hbox{%
$\scriptstyle #1$}\kern \wd1{\mathrm{C}}_{\kern -1.05\wd0\kern -0.99\wd1{#1}
 \kern 1.15\wd0{#2}}}

\def\clvec#1#2#3{\def\clvecone{#3}\left(\arrb{c} \dsp #1\\ \dsp #2
 \ifx\clvecone\empty\else\\ \dsp #3\fi\arre\right)}

\def\pderiv#1#2{\dsp\frac{\partial\,#1}{\partial#2}}

\def\le{\leqq} \def\leq{\leqq}\def\ge{\geqq} \def\geq{\geqq}

\def\reals{{\mathbb R}}

\def\preals{[0,\infty)} 
\def\pintegers{{\mathbb Z}_+}

\def\prb#1{\def\prbone{#1}
  \ifx\prbone\empty{\mathrm{P}}\else{\mathrm{P[\;}}#1{\mathrm{\;]}}\fi}
\def\prbseq#1#2{\def\prbseqone{#2}
  \ifx\prbseqone\empty{\mathrm{P}}_{#1}\ignorespaces
  \else{\mathrm{P}}_{#1}{\mathrm{[\;}}#2{\mathrm{\;]}}\fi}

\def\EEseq#1#2{\def\EEseqone{#2}
  \ifx\EEseqone\empty{\mathrm{E}}_{#1}\else
 {\mathrm{E}}_{#1}{\dsp\mathrm{[\;}}#2{\mathrm{\;]}}\fi}

\def\VVseq#1#2{\def\VVseqone{#2}
  \ifx\VVseqone\empty{\matrm{V}}_{#1}\else
 {\mathrm{V}}_{#1}{\dsp\mathrm{[\;}}#2{\mathrm{\;]}}\fi}

\def\ssN{^{(N)}}

\def\chrfcn#1{\mathop{\mathbf{1}}\nolimits_{#1}}

\title{
 Stochastic ranking process with space-time dependent intensities
}
\author{
Tetsuya Hattori
\footnote{ 
Partly supported by KAKENHI 21340020 (the Grant-in-Aid for 
Scientific Research (B)) from the Japan Society for the Promotion of Science. 
}
\\ 
\small Laboratory of Mathematics, Faculty of Economics, Keio University, 
\\
\small Hiyoshi Campus, 4--1--1 Hiyoshi, Yokohama 223-8521, Japan
\\ \small URL: \url{http://web.econ.keio.ac.jp/staff/hattori/research.htm}
\\ \small email: \url{hattori@econ.keio.ac.jp}
\\ \and
Seiichiro Kusuoka
\footnote{ 
Research Fellow of the Japan Society for the Promotion of Science 
}
\\
\small 
Mathematical Institute, Graduate School of Science, Kyoto University,
\\
\small Kita-Shirakawa, Sakyo-ku, Kyoto 606-8502, Japan
\\ \small email: \url{kusuoka@math.kyoto-u.ac.jp}
} 
\date{2011/11/24}

\begin{document}

\maketitle

\begin{abstract}
We consider the stochastic ranking process
with space-time dependent jump rates for the particles.
The process is a simplified model of the time evolution of the rankings 
such as sales ranks at online bookstores.
We prove that the joint empirical distribution of jump rate and 
scaled position converges almost surely to a deterministic distribution, 
and also the tagged particle processes converge almost surely,
in the infinite particle limit.
The limit distribution is characterized by a 
system of inviscid Burgers-like integral-partial differential equations
with evaporation terms,
and the limit process of a tagged particle is a motion along a characteristic
curve of the differential equations except at its Poisson times of jumps to 
the origin.
\end{abstract}

\noindent 
2000 \textit{Mathematics Subject Classification}.
Primary 60K35;
Secondary 35C05,
82C22.

\noindent 
\textit{Key words}. 
Stochastic ranking process, Poisson process, hydrodynamic limit,
inviscid Burgers equation, move-to-front rules.

\newpage

\setcounter{section}{0}

\setcounter{equation}{0}

\section{Introduction.}
\seca{1}
In this paper, we consider stochastic ranking processes whose jump rates 
depend not only on time but also on their positions.
Stochastic ranking processes are a model of a ranking system,
such as the sales ranks found at online bookstores.
We consider $N$ particles each of which are exclusively located 
at $1,2,\ldots ,N$.
Each particle jumps to $1$ according to its Poisson clock.
When a jump of the particle at position $i$ occurs, the particle 
moves to position $1$ and the locations of the particles at $1,2,\dots ,i-1$
are sifted by $+1$.
Particles whose Poisson clocks rang recently are 
at positions with small numbers,
and the others are at positions with large numbers.
We regard the number for each particle as the particle's rank.
This system enables us to give ranks to $N$ particles, 
and in this paper we call the time evolution of the particles 
given by this ranking system the \textit{stochastic ranking process}.

A precise formulation of the stochastic ranking process 
which we consider in this paper are as follows.
Let $(\Omega,{\cal F},P)$ be a probability space,
and let $\{ \nu _i (d\xi ds)\} _{i=1,2,3,\dots }$ be independent 
Poisson random measures on $[0,\infty )\times[0,\infty )$ 
with the intensity measure $d\xi ds$.
Let $W$ be a set of non-negative valued $C^1$ functions
\[ w:\ [0,1]\times\preals\to\preals, \]
such that, for each $T>0$,
\eqnb
\eqna{rw}
\eqna{yCy0eq0prf3}
R_w(T):=\sup_{w\in W} \sup_{(y,t)\in[0,1]\times[0,T]} 
\max\biggl\{
w(y,t),\ \left|\pderiv{w}{y}(y,t)\right| \biggr\} < \infty.
\eqne

Let $w_i$, $i=1,2,\ldots$ be a sequence in $W$, and for 
a positive integer $N$, put
\[
w\ssN_{i}(k,t):= w_{i}(\frac{k-1}{N},t),\ \ 
k=1,2,\ldots,N,\ t\in [0,\infty),\ i=1,2,\ldots,N.
\]
Also, let $x^{(N)}_1$, $x^{(N)}_2$, $\ldots$, $x^{(N)}_N$ be 
a rearrangement of $1,2,\dots,N$.
Define a process 
\[X\ssN=(X\ssN_1,\ldots,X\ssN_N)\]
by
\eqnb
\eqna{SIformSRP}
\arrb{l}\dsp
X\ssN_i(t) \\
\dsp = x\ssN_i + \sum _{j=1} ^N
 \int _{s\in (0,t]} \int _{\xi \in [0,\infty )}
 \chrfcn{ X\ssN_i(s-)<X\ssN_j(s-) }
 \chrfcn{\xi\in [0,w\ssN_{j} (X\ssN_j(s-),s))}
 \nu _j(d\xi ds)
\\ \dsp \phantom{=}
 + \int _{s\in (0,t]}  \int _{\xi \in [0,\infty )} (1-X\ssN_i(s-))
\chrfcn{\xi\in [0,w\ssN_{i}(X\ssN_i(s-),s))} 
 \nu _i(d\xi ds), \ \ i=1,2,\dots ,N,\ t\geq 0,
\arre
\eqne
where, $\chrfcn{B}$ is the indicator function of event $B$.
The integrands in the \eqnu{SIformSRP} are predictable,
hence the right hand side of \eqnu{SIformSRP} is well-defined as
the Ito--integrals \cite[\S \Rom{4}.9]{IW}.

$X\ssN(t)$ is a rearrangement of $1,2,\ldots,N$ for all $t\ge 0$,
which we regard as ranks or positions of particles $1,2,\ldots,N$
at time $t$.
Moreover, for $i=1,2,\ldots,N$, and $t>t_0\ge 0$, let
\eqnb
\eqna{Jtopnu}
\eqna{Jtopist}
J\ssN_i(t_0,t)=\left\{  
\int_{s\in (t_0,t]} \int _{\xi \in [0,\infty )}
 \chrfcn{\xi\in [0,w\ssN_{i}(X\ssN_j(s-),s))} \nu _i(d\xi ds) > 0
\right\}.
\eqne
Then, the last term on the right hand side of \eqnu{SIformSRP}
implies that $J\ssN_i(t_0,t)$ denotes the event
that the particle $i$ jumps to the top position ($X\ssN_i(s)=1$)
in the time interval $(t_0,t]$.
Also, the second term on the right hand side of \eqnu{SIformSRP}
implies that conditioned on the complement $J\ssN_i(t_0,t)^c$ of
$J\ssN_i(t_0,t)$,
\eqnb
\eqna{Jstepdown}
X\ssN_i(s)-X\ssN_i(s-)=0\ \mbox{ or }\ 1,
\eqne
for $t_0<s\le t$, where the latter occurs if and only if a particle $k$ 
at tail side ($X\ssN_k(s)>X\ssN_i(s)$) jumps to top at time $s$.

With a great advance in the internet technologies,
a new application of the process appeared.
The ranking numbers such as those found in the web pages of online retails, 
e.g., the sales ranks of books at the Amazon online bookstore,
are found to follow the predictions of 
the model \cite{HH072,mv2frnt,ranking,dojin}.
In the ranking of books,
each time a book is sold its ranking spontaneously jumps to small numbers
(relatively close to $1$), regardless of how bad its previous position
was (large $X\ssN_i(t-)$ in our notation), and
regardless of how unpopular (small $w\ssN_i$, in our notation) the book is.
The stochastic process we consider here corresponds to a mathematical 
simplification of this observation, that each time a book is sold
its ranking jumps to $1$ instantaneously.
With a view that the process is a model of such on-line, real time, 
rankings of a large number of items according to their popularity,
we will call the model the stochastic ranking processes.

At first thought one might guess that such a naive ranking rules
of spontaneous jump to $1$ at each sale, as in the definition of
the stochastic ranking processes,
will not be a good index for the popularity of books.
But with a closer look, one notices that the well sold books are dominant
near the top position, while books near the tail position are rarely sold.
Though the rankings of each book are stochastic with sudden jumps,
the spatial distribution of the jump rates is more stable.
In the bookstore's view, what matters is not a specific book,
but the totality of sales.
This motivates an interest on the evolution of the 
joint empirical distribution of position and jump rates.

In \cite{HH071,HH072,timedep,Nagahata10}, infinite particle
(large $N$) scaling limit for this model is considered, 
and the explicit formula of the limit distribution of the joint
empirical distribution of scaled position and the jump rate is found, 
which further is characterized as a solution to a system of 
inviscid Burgers-like equations with a term representing evaporation.
The limit formula is successfully applied to the time developments
of ranking numbers such as those found in the web pages of online 
bookstores \cite{HH072,mv2frnt,ranking}.
Furthermore, convergence of the joint empirical distribution as
a process and convergence of tagged particle process
are proved in \cite{Nagahata10}.

If the model \eqnu{SIformSRP} is independent of spatial position, 
i.e., if $w\ssN$'s are independent of their first variables $x$,
then the law of the process \eqnu{SIformSRP} reduces to 
that of \cite[eq.~(2)]{timedep} and \cite[eq.~(1)]{Nagahata10},
the stochastic ranking process with time dependent (but position independent)
intensities.
Thus \eqnu{SIformSRP} is an extension of \cite{timedep,Nagahata10} to the case
where the dynamics is dependent on the value of $X\ssN_i(t)$\,,
i.e., to the position dependent case.
In the present paper, 
we mathematically extend the previous results to the case where the jump rates
$w\ssN$'s are both position and time dependent.

If $w\ssN$'s are positive constants,
\eqnu{SIformSRP} further reduces to the homogeneous case considered
in \cite{HH071,HH072}.
A discrete time version of the homogeneous case
has been known since \cite{Tsetlin1963},
and has been extensively studied since then
and is called move-to-front (MTF) rules
\cite{mv2frnt1,mv2frnt2,mv2frnt3,Letac1974,Kingman75}.
The process and its generalization have, in particular, 
been extensively studied as a model of 
least-recently-used (LRU) caching in the field of information theory 
\cite{Rivest1976,Fagin77,Bitner79,CHS88,BlomHolst91,Rodrigues95a,Fill96JTP,%
mv2frnt4,Fill96TCS,Jelenkovic99,SM2006,Jelenkovic08,BF10},
and also is noted as a time-reversed process of top-to-random shuffling.

A motivation for an on-line web retail store to provide the sales ranks,
in their web pages for public access,
would be to give information on the popularity of each products
which the store provides, to attract consumers' attention on popular
products. 
Extending previous results to the case of position dependent jump rates,
which is the main aim of the present paper,
corresponds to providing a mathematical framework for considering a possibility
of such expected effect of popular products receiving extra attention
and effectively increase their jump rates according to their rankings.

We introduce the normalized position for each particle $i$ at time $t$
\eqnb
\eqna{YN}
 Y\ssN_i(t)=\frac1N(X\ssN_i(t)-1),
\end{equation}
and consider the joint empirical distribution of jump rate
and normalized position, given by
\eqnb
\eqna{mut}
\mu\ssN_{t} =\frac1N \sum_{i=1}^N \delta_{(w_i,Y\ssN_i(t))}\,.
\eqne
(We will denote a unit measure on any space by $\delta_c$.)
$\mu\ssN_{\cdot}$ is a stochastic process
taking values in the set of Borel probability measures.
For each $T>0$, $\mu\ssN_t$, $t\in[0,T]$, is regarded as a 
stochastic process on $C^{1,0}([0,1]\times[0,T])\times[0,1]$, 
where $C^{1,0}([0,1]\times[0,T])$ is the total set of functions 
$f\in C([0,1]\times[0,T])$ 
such that $\pderiv{f}{y}(y,t)\in C([0,1]\times[0,T])$.
$C^{1,0}([0,1]\times[0,T])$ is a Polish space (complete separable
metric space) with norm
\[ \sup_{(y,t)\in[0,1]\times[0,T]} \biggl\{  |w(y,t)|,\ 
\left|\pderiv{w}{y}(y,t)\right|\ \biggr\}. \]
Since $C^1([0,1]\times[0,T])$  is a Polish space, 
so is $C^1([0,1]\times[0,T]) \times[0,1]$ \cite[Example 26.2]{Bauer}.
We assume a standard topology of weak convergence of probability measures
on $C^1([0,1]\times[0,T])\times[0,1]$.

To prove convergence of measures, we work with a distribution function.
For each integer $N$ define
\eqnb
\eqna{distributionfcn}
\eqna{distributionfcnexplicit}
\ \ \ 
U\ssN(dw,y,t)=\mu\ssN_{t}(dw \times [y,1))
=\frac1N \sum_{i=1}^N \chrfcn{X\ssN_i(t-)\ge Ny+1}\,\delta_{w_i}(dw)\,,
\\ 0\le y\le 1,\ t\ge0.
\eqne
For each $(y,t)$, $U\ssN(\cdot,y,t)$ is a Borel measure on $W$.
Note that $U\ssN(dw,y,t)$ is non-increasing in $y$ and satisfies
\eqnb
\eqna{discreteU}
\int_W U\ssN(dw,y,t) = \frac{[N\,(1-y)]}N\,,
\ \ 0\le y\le 1,\ t\ge0.
\eqne
where, for real $z$, $[z]$ is the largest integer not exceding $z$.

As an extension of the corresponding results in \cite{HH072,timedep}, 
the infinite particle scaling limit $U$ of
$U\ssN$ turns out to be characterized by a 
system of inviscid Burgers-like integral--partial differential equations
with evaporation terms.
Denote the set of `boundary points' and of `initial points' by
\eqnb
\eqna{bdrypts}
\Gamma_b=\{(0,t_0)\mid t_0\ge0\},
\eqne
and
\eqnb
\eqna{initpts}
\Gamma_i=\{(y_0,0)\mid 0\le y_0\le 1\},
\eqne
respectively, and put
\eqnb
\eqna{spacetimeboundary}
\Gamma=\Gamma_b\cup\Gamma_i\,.
\eqne
Also, for $t\ge 0$ put
\eqnb
\eqna{Gammat}
 \Gamma_t=\{(y_0,t_0)\in \Gamma \mid t_0\le t\}
=\Gamma_i\cup \{(0,t_0) \mid 0\le t_0\le t\}.
\eqne
\thmb
\thma{Burgers}
Let $\lambda$ be a Borel probability measure on $W$, and
$\rho:\ W\times[0,1]\to [0,1]$ be a non-negative Borel measurable
function such that $\dsp\pderiv{\rho}{y}(w,y)$ exists and continuous, and
\eqnb
\eqna{rhoasdistributionfunction}
\pderiv{\rho}{y}(w,y)\le 0,\ (w,y)\in W\times[0,1],
\eqne
and such that $\rho(w,0)=1$ and $\rho(w,1)=0$, $w\in W$.
Define a Borel measure on $W$ parametrized by $y\in[0,1]$, by
\eqnb
\eqna{U0yabscontiU00}
U_0(dw,y)=\rho(w,y)\,\lambda(dw),\ \ y\in[0,1],\ w\in W.
\eqne
In particular, $U_0(dw,0)=\lambda(dw)$.
Assume also
\eqnb
\eqna{Burgersassump}
 U_0(W,y)=\int_W U_0(dw,y) = 1-y,\ 0\le y\le 1.
\eqne

Then there exists a unique pair of functions
\[
y_C:\ \{(\gamma,t)\in\Gamma\times \preals\mid \gamma\in\Gamma_t\}\to[0,1],
\]
and $U=U(dw,y,t)$ on $[0,1]\times\preals$ taking values in
the non-negative Borel measures on $W$,
such that, 
\itmb
\item
$y_C(\gamma,t)$ and $\dsp \pderiv{y_C}{t}(\gamma,t)$ is continuous
\item
for each $t>0$,
$\dsp y_C(\cdot,t):\ \Gamma_t\to[0,1]$ is surjective,
\item
for all bounded continuous $h:\ W\to\reals$,
$\dsp U(h,y,t):=\int_W h(w) U(dw,y,t)$ is 
Lipschitz continuous in $(y,t)\in [0,1]\times[0,T]$ for any $T>0$,
and non-increasing in $y$, and
\item
the following \eqnu{U0}, \eqnu{phic}, \eqnu{yc}, and \eqnu{soliditycond}
hold:
\itme
\eqnb
\eqna{U0}
y_C(\gamma,t_0)=y_0
\ \mbox{ and }\ 
 U(dw,y_0,t_0)=U_0(dw,y_0),\ \ \gamma=(y_0,t_0)\in\Gamma,
\eqne
\eqnb
\eqna{phic}
U(h,y_C(\gamma,t),t)=U_0(h,y_0)-\int_{t_0}^t V(h,y_C(\gamma,s),s)\,ds,
\ \ t\ge t_0,\ \gamma=(y_0,t_0)\in \Gamma,
\eqne
for all bounded continuous function $h:\ W\to\reals$,
where, $\dsp U(h,y,t):=\int_W h(w) U(dw,y,t)$, and
\eqnb
\eqna{prfofHDL3}
V(h,y,t)=\int_W h(w)\,w(y,t)\,U(dw,y,t)
+\int_y^1 \int_W h(w)\, \pderiv{w}{z}(z,t)\,U(dw,z,t)\,dz,
\eqne
and
\eqnb
\eqna{yc}
\pderiv{y_C}{t}(\gamma,t)= V(\chrfcn{W}, y_C(\gamma,t),t),
\ \ t\ge t_0,\ \gamma=(y_0,t_0)\in \Gamma,
\eqne
where 
$\chrfcn{W}(w)=1$ for all $w\in W$,
and
\eqnb
\eqna{soliditycond}
U(\chrfcn{W},y,t)=1-y,\ \ 0\le y\le 1,\ t\ge0.
\eqne
\DDD
\thme

The claim \eqnu{soliditycond}, together with continuity and monotonicity 
of $U$, implies that $U$ determines a Borel probability measure $\mu_t$ 
on the direct product $W\times [0,1]$  parametrized by $t$:
\eqnb
\eqna{distributionfcncontilimit}
U(dw,y,t)=\mu_{t}(dw \times [y,1)),\ \ 0\le y\le 1,\ t\ge0.
\eqne

If $U(h,y,t)$ in \thmu{Burgers} is $C^1$ in a neighborhood of 
$(y,t)\in (0,1)\times (0,\infty)$, 
then differentiating \eqnu{phic} by $t$ and using \eqnu{yc},
and noting that $y_C(\cdot,t):\ \Gamma_t\to [0,1]$ is surjective,
we have 
\eqnb
\eqna{Burgers}
\pderiv{U}{t}(h,y,t) 
+ V(\chrfcn{W},y,t)\, \pderiv{U}{y}(h,y,t)
=-V(h,y,t),
\eqne
where $V$ is as in \eqnu{prfofHDL3}.
$y_C$ in \eqnu{yc} determines the characteristic curves for \eqnu{Burgers}.
In terms of \cite[\S 3.4]{Bressan},
we can therefore say that \thmu{Burgers} claims 
global existence of the Lipschitz solution
(broad solution which is Lipschitz continuous) to the system of
quasilinear partial differential equations \eqnu{Burgers},
with components parametrized by (possibly continuous) parameter $w$.
To be more precise, 
we have extended the definition in \cite[\S 3.4]{Bressan}
of Lipschitz solution for \eqnu{Burgers}
to the non-local case (see \eqnu{prfofHDL3}),
and for the case where $V(\chrfcn{W},y,t)$ in the
left-hand side of \thmu{Burgers} is common for all $h$.
We have also generalized the notion of domain of determinancy
defined in \cite[\S 3.4]{Bressan}, which in the present case corresponds to 
\[ \{(y,t)\in[0,1]\times\preals\mid y\ge y_C((0,0),t)\}, \]
to the domain determined by boundary conditions
\[ \{(y,t)\in[0,1]\times\preals\mid y<y_C((0,0),t)\}, \]
with initial data $U(h,\cdot,0)=U_0(h,\cdot)$ and 
the boundary condition $U(h,0,t)=U_0(h,0)$, $t\ge0$,
as obtained in \eqnu{U0}.

As a simple example, where the jump rates are finitely many
space-time constants, there is a natural one to one onto
map from $W$ (the space of jump rate functions), to 
a finite set $\{w_1,w_2,\ldots,w_A\}$ of positive integers
for some positive integer $A$,
and the distribution of the jump rates $U_0(dw,0)=\lambda(dw)$
can be identified with 
\[
\lambda(dw)=\sum_{a=1}^A r_a\,\delta_{w_a},
\]
for some positive constants $r_a$, $a=1,2,\ldots,A$,
satisfying $\dsp \sum_{a=1}^A r_a=1$.
In this simple example, \eqnu{Burgers} reduces to
\eqnb
\eqna{Burgersfinitetypes}
\pderiv{U_a}{t}(y,t) 
+ \sum_{b=1}^A w_b U_b(y,t)\, \pderiv{U_a}{y}(y,t)
=-w_a U_a(y,t),
\eqne
where we wrote $U_a(y,t)=U(h,y,t)$ and
$V_a(y,t)=V(h,y,t)=w_a U_a(y,t)$,
so that
$\dsp V(\chrfcn{W},y,t)=\sum_{b=1}^A w_b U_b(y,t)$.
If the right hand side of \eqnu{Burgersfinitetypes} is $0$,
the partial differential equation is known (for $A=1$) as
the inviscid Burgers equation in the terminology of fluid dynamics.
In terms of fluid dynamics, the right hand side of 
\eqnu{Burgersfinitetypes} could be interpreted as the evaporation
of the fluid.

For the case \eqnu{Burgers} which we consider in this paper,
non-locality of interaction is inevitable, precisely because of
the position dependence of the jump rate functions,
hence we need to consider a harder problem of a system of 
differential--integral equations compared to previous cases 
\cite{HH071,timedep}.

Now we give a norm of measures in order to state the next theorem.
Let $||\cdot ||_{\rm{var}}$ be the total variation norm for Borel measures on 
$W$, i.e. for a signed measure $\mu$ on $W$ define $||\mu ||_{\rm{var}}$ by
\[
||\mu ||_{\rm{var}} = \mu ^+(W) + \mu ^-(W),
\]
where $\mu ^+$ and $\mu ^-$ are the positive part and the negative part 
obtained by Hahn-Jordan decomposition of $\mu$ respectively.

We consider a scaling limit of the stochastic ranking process as $N\to\infty$ (the limit for the number of particles to infinity).
We are naturally considering
a law of large number type of results, and as suggested by the fact,
which we state in the following, that the limit distribution satisfies 
\eqnu{phic}, or more intuitively, non-linear equations \eqnu{Burgers}, 
it is a non-trivial problem of 
the law of large numbers for dependent variables.
This has been the case also for previous results in \cite{HH071,timedep},
but in the previous studies, where the jump rate functions are
independent of spatial positions,
a special combination of quantities we define
($U^{(N)}(B,Y^{(N)}_C(y_0,t_0,t),t)$, in terms of notations in
\secu{profofHDL}) turns out to be a sum of independent random variables.

However, the position-dependence of jump rates, as considered in
the present paper, implies that the dependence of random variables are
built-in in the model, so that the proofs in \cite{HH071,timedep} do not work
in the present case.
Inspired partly by \cite{Nagahata10}, where the case of 
finite types of position independent particles are proved 
\cite[Prop.~1.1 and Thm.~1.2]{Nagahata10},
we extend his result to our position-dependent case, and obtain a convergence of empirical distribution and also
the limiting dynamics of fixed finite particles (tagged particles)
for the case of jump rate functions with space-time dependence as follows.
\thmb
\thma{HDLposdep}
Assume that with probability $1$, 
\eqnb
\eqna{Uinitconvass}
 \limf{N} \sup _{y\in[0,1)}||U\ssN(\cdot ,y,0) - U_0(\cdot ,y)||_{\rm{var}}=0,
\eqne
where $U_0(dw,y)$ satisfies all the assumptions in \thmu{Burgers}.
Then the following hold.
\itmb
\item
With probability $1$, for all $T>0$,
$\dsp \limf{N} U\ssN(dw,y,t)= U(dw,y,t)$, 
uniformly in $y\in[0,1)$ and $t\in[0,T]$,
where $U$ is the solution claimed in \thmu{Burgers}.

\item
Assume in addition that, 
\eqnb
\eqna{taggptclconvass}
\limf{N} \frac1N x\ssN_{i} =y_i\,,\ \ i=1,2,\ldots,L,
\eqne
for a positive integer $L$ and $y_i\in[0,1)$, $i=1,2,\ldots,L$.
Then, with probability $1$, for all $T>0$,
the tagged particle system 
\[
(Y\ssN_1(t),Y\ssN_2(t),\ldots,Y\ssN_L(t))
\]
converges as $N\to\infty$, uniformly in $t\in[0,T]$ to a limit
$(Y_1(t),Y_2(t),\ldots,Y_L(t))$. Here, for each $i=1,2,\ldots,L$, 
$Y_i$ is the unique solution to
\eqnb
\eqna{taggedlimit}
\arrb{l}\dsp
Y_i(t)
\\\dsp{}
 = y_i + \int _0^t V(\chrfcn{W},Y_i (s-),s) ds
 - \int _{s\in (0,t]} \int _{\xi \in [0,\infty )} Y_i(s-) 
\chrfcn{\xi \in [0,w_i (Y_i(s-),s))} \nu _{i}(d\xi ds),
\arre
\eqne
where, $V$ is as in \eqnu{prfofHDL3}.
\itme
\DDD\thme

When $\{ w_i; i=1,2,3,\dots \}$ is a finite set of $W$, because of 
Proposition \ref{appendix} in Appendix, we obtain the following corollary 
easily.

\corb
When $w_i \in \{ \tilde w_\alpha \in W; \alpha =1,2,\dots ,A\}$ 
for $i=1,2,3,\dots$,
the assumption \eqnu{Uinitconvass} of \thmu{HDLposdep} is relaxed 
as follows:
\[
 \limf{N} U\ssN(\{ \tilde w_\alpha\},y,0)
= U(\{ \tilde w_\alpha\} ,y,0),\ \mbox{ for each }\ y\in[0,1)
\]
with probability $1$ for $\alpha =1,2,\dots ,A$.
\core

A discrete correspondence $Y\ssN_C$ of the characteristic curves
$y_C$ is defined in \eqnu{yCN} in \secu{profofHDL},
both of which have been key quantities since \cite{HH071}.
As pointed in \cite{Nagahata10}, 
\eqnu{taggedlimit} says that a particle moves along a characteristic curve
of \eqnu{Burgers} except at its Poisson times of jumps to $y=0$.

\bigskip\par

The plan of the paper is as follows.
In \secu{prfofBurgers} we prove \thmu{Burgers},
and in \secu{profofHDL} we prove \thmu{HDLposdep}.

\bigskip\par\textbf{Acknowledgment.}
The authors would like to thank Prof.~Y.~Nagahata for discussions.
T.H. also would like to thank Prof.~M.~Hino, Prof.~I.~Shigekawa, 
Prof.~S.~Takesue,
Prof.~K.~Yano, Prof.~Y.~Yano, and Prof.~N.~Yoshida, for their interest
and discussions on the present work, and also for their hospitality
at Kyoto University.

\section{Proof of \protect\thmu{Burgers}.}
\seca{prfofBurgers}
Consider first the case $(y_0,t_0)\in\Gamma_i$, namely, the case $t_0=0$.
\lemb
\lema{yCt0eq0}
There exists a unique $C^1$ function $f:\ [0,1]\times\preals\to[0,1]$
which satisfies
\eqnb
\eqna{yCt0eq0}
f(y,t)
=1+\int_W \biggl(\int_y^1 \pderiv{\rho}{z}(w,z)
 \,\exp(-\int_0^t w(f(z,s),s)\,ds)\,dz \biggr) U_0(dw,0),
\ y\in[0,1],\ t\ge0,
\eqne
where $\rho$ and $U_0$ are as in the assumptions of \thmu{Burgers}.
\DDD\leme
\prfb
For $k\in\pintegers$, define $f_k:\ [0,1]\times \preals\to[0,1]$ 
inductively by
\[
f_0(y,t)=1+\int_W \biggl(\int_y^1 \pderiv{\rho}{z}(w,z)
 \,\exp(-\int_0^t w(z,s)\,ds)\,dz \biggr)\,U_0(dw,0),
\]
and
\eqnb
\eqna{yCt0eq0prf2}
f_{k+1}(y,t)
=1+\int_W \biggl(\int_y^1 \pderiv{\rho}{z}(w,z)
 \,\exp(-\int_0^t w(f_k(z,s),s)\,ds)\,dz \biggr)\,U_0(dw,0),
\ \ k\in\pintegers.
\eqne
Assume that $f_k$ is continuous and takes values in $[0,1]$.
Then \eqnu{yCt0eq0prf2} is well-defined.
Non-increasing assumption of \thmu{Burgers} for $\rho$ implies
$\dsp \pderiv{\rho}{z}(w,z)\le 0$, 
hence \eqnu{yCt0eq0prf2} implies $f_{k+1}\le 1$.
Similarly, using also \eqnu{U0yabscontiU00} and \eqnu{Burgersassump},
\[
\arrb{l}\dsp
f_{k+1}(y,t)
\\ \dsp {}
\ge 1+\int_W \biggl(\int_y^1 \pderiv{\rho}{z}(w,z)\,dz \biggr)\,U_0(dw,0)
= 1+\int_W U_0(dw,1)-\int_W U_0(dw,y) = y
\\ \dsp {}
\ge 0.
\arre
\]
$\rho$ and $w$ are $C^1$ in $z$, by assumption of \thmu{Burgers},
hence \eqnu{yCt0eq0prf2} implies that $f_{k+1}$ is continuous.
By induction, $f_k$ is continuous and takes values in $[0,1]$, for all $k$.

For $k\in\pintegers$, put $F_k(y,t)=|f_{k+1}(y,t)-f_{k}(y,t)|$. 
Then, using \eqnu{rw} and the assumptions of \thmu{Burgers} as above,
we have
\eqnb
\eqna{yCt0eq0prf3}
F_{k+1}(y,t)\le R_w(T) \int_y^1 \int_0^t F_{k}(z,s)\,ds\, dz,
\ \ y\in[0,1],\ t\in[0,T],\ k\in\pintegers,
\eqne
for any $T>0$.
Since all $f_k$'s are continuous and take values in $[0,1]$,
$F_{k}$, $k=1,2,\ldots$, are also continuous and take values in $[0,1]$.
Then it holds by the argument of \cite[\S 3.8, Lemma 3.4]{Bressan}, that
\eqnb
\eqna{iterationgood}
0\le F_{k}(y,t)\le e^{2R_w(T)t}2^{-k}
,\ \ y\in[0,1],\ t\in[0,T],\ k\in\pintegers.
\eqne
In fact, since $F_0$ takes values in $[0,1]$, 
\eqnu{iterationgood} holds for $k=0$.
Assume \eqnu{iterationgood} holds for some $k$.
Then \eqnu{yCt0eq0prf3} implies
\[
F_{k+1}(y,t)\le 2^{-k} R_w(T)\int_y^1 \int_0^t e^{2R_w(T)s}\,ds
\le e^{2R_w(T)t}2^{-k-1},\ \ 0\le y<1,\ 0\in[0,T].
\]
By induction, \eqnu{iterationgood} holds for all $k\in\pintegers$.
In particular,
$\dsp f_0(y,t)+\sum_{k=0}^{\infty} F_{k}(y,t)$
converges uniformly in $(y,t)$ for any bounded range of $t$.
Hence,
$\dsp f_k(y,t)=f_0(y,t)+\sum_{j=0}^{k-1} (f_{j+1}(y,t)-f_j(y,t))$
converges as $k\to \infty$ to a function, continuous in $y$ and $t$.
Let
\[ f(y,t)=\limf{k} f_k(y,t),\ \ y\in[0,1],\ t\ge0. \]
Then \eqnu{yCt0eq0prf2} implies that $f$ satisfies \eqnu{yCt0eq0}.
Also, $0\le f_k\le 1$ implies
\eqnb
\eqna{frange}
0\le f(y,t)\le 1,\ 0\le y\le 1,\ t\ge 0.
\eqne
The right hand side of \eqnu{yCt0eq0}, with the assumptions in
\thmu{Burgers} implies that $f(y,t)$ is $C^1$.

Next, we prove the uniqueness.
Suppose for $i=1,2$,
$f^{(i)}:\ [0,1]\times\preals\to[0,1]$ are continuous functions
which satisfy \eqnu{yCt0eq0}. Then 
$|f^{(1)}(y,0)-f^{(2)}(y,0)|=0$ and, as above, for each $T>0$,
\[
|f^{(1)}(y,t)-f^{(2)}(y,t)|
\le R_w(T) \int_y^1\int_0^t|f^{(1)}(z,s)-f^{(2)}(z,s)|\,ds\,dz
\ \ y\in[0,1],\ t\in[0,T], \]
which implies $f^{(1)}=f^{(2)}$.
\QED\prfe

Next, consider the case $(y_0,t_0)\in\Gamma_b$, namely, the case $y_0=0$.
\lemb
\lema{yCy0eq0prf1}
For each continuous function 
$\tilde{g}: \{(s,t)\in\preals^2\mid 0\le s\le t\}\to [0,1]$,
there exists a unique non-negative function
$\dsp \eta:\ W\times \preals\to\preals$, integrable with respect to 
$U_0(dw,0)$, continuous in the second variable, 
which satisfy, for each $w\in W$,
\eqnb
\eqna{yCy0eq0prf1}
\arrb{l}\dsp
\eta(w,t)
=\int_0^t \eta(w,u)\, w(\tilde{g}(u,t),t)
\, \exp(-\int_u^t w(\tilde{g}(u,v),v)\,dv) \, du
\\ \dsp \phantom{\eta(w,t)=}
- \int_0^1 \pderiv{\rho}{z}(w,z) \, w(f(z,t),t)
\, \exp(-\int_0^t w(f(z,v),v)\,dv) \, dz, \ \ t\ge0,
\arre
\eqne
where $\rho$ is as in the assumption of \thmu{Burgers}, 
and $f$ is the function given by \lemu{yCt0eq0}.

Moreover, it holds that
\eqnb
\eqna{yCy0eq0prf2}
\arrb{l}\dsp
\int_0^t \eta(w,u)
\, \exp(-\int_u^t w(\tilde{g}(u,v),v)\,dv) \, du
\\ \dsp {}
=1
+\int_{0}^1 \pderiv{\rho}{z}(w,z) 
\exp(-\int_0^t w(f(z,s),s)\,ds) \,dz.
\arre
\eqne
In particular,
for any $T>0$, there exists $C(T)>0$, which is independent of
$\tilde{g}$, such that 
\eqnb
\eqna{yCy0eq0prf2b}
0\le \int_W \eta(w,t)\, U_0(dw,0)\le C(T),\ 0\le t\le T.
\eqne
\DDD\leme
\prfb
Define a sequence of continuous functions
$\dsp \eta_{k}:\ W\times\preals\to\preals$, $k=0,1,2,\ldots$,
inductively, by
\[ \eta_{0}(w,t)=0,\ \ w\in W,\ t\ge 0,\]
and
\eqnb
\eqna{etak}
\arrb{l}\dsp
\eta_{k+1}(w,t)
=\int_0^t \eta_{k}(w,u)\, w(\tilde{g}(u,t),t)
\, \exp(-\int_u^t w(\tilde{g}(u,v),v)\,dv) \, du
\\ \dsp \phantom{\eta_{k+1}(w,t)=}
- \int_0^1 \pderiv{\rho}{z}(w,z) \, w(f(z,t),t)
\, \exp(-\int_0^t w(f(z,v),v)\,dv) \, dz.
\arre
\eqne
For $k\in\pintegers$
put $\dsp H_{k}(t)=\int_W |\eta_{k+1}(w,t)-\eta_{k}(w,t)|\, U_0(dw,0)$.
Non-negativity of $w\in W$ and \eqnu{yCy0eq0prf3} imply
\[ H_{k+1}(t)\le R_w(T)\int_0^t H_{k}(u)\,du,\ 0\le t\le T. \]
\cite[\S 3.8, Lemma 3.4]{Bressan} implies that there exists a positive
constant $C(T)$ such that
\[ H_k(t) \le C(T)\, 2^{-k},\ t\in[0,T],\ k\in\pintegers, \]
hence, as in the proof of \lemu{yCt0eq0},
$\dsp \eta(w,t)=\limf{k} \eta_{k}(w,t)$ exists,
is continuous, non-negative, and satisfies \eqnu{yCy0eq0prf1}.
Integrability inductively follows from \eqnu{etak} by
\[ \arrb{l}\dsp
\sup_{t\in[0,T]}\int_W \eta_{k+1}(w,t)\,U_0(dw,0)
\\ \dsp {}
\le R_w(T)\,\sup_{t\in[0,T]}\int_W \int_0^t \eta_{k}(w,u) \, du\,U_0(dw,0)
+ R_w(T)\,\int_W 
\int_0^1 \biggl(-\pderiv{\rho}{z}(w,z)\biggr) \, dz\,U_0(dw,0)
\\ \dsp {}
=R_w(T)\,\sup_{t\in[0,T]}\int_W \int_0^t \eta_{k}(w,u) \, du\,U_0(dw,0)
+ R_w(T),
\arre \]
where we also used 
\eqnu{rhoasdistributionfunction},
\eqnu{U0yabscontiU00} and \eqnu{Burgersassump}.

Next, we prove the uniqueness.
Suppose for $i=1,2$,
$\eta^{(i)}:\ W\times\preals\to\preals$ are functions,
continuous in the second variable and satisfy \eqnu{yCy0eq0prf1}. Then 
$|\eta^{(1)}(w,0)-\eta^{(2)}(w,0)|=0$
and, as above, for each $T>0$,
\[
|\eta^{(1)}(w,t)-\eta^{(2)}(w,t)|
\le R_w(T) \int_0^t|\eta^{(1)}(w,s)-\eta^{(2)}(w,s)|\,ds
\ \ t\in[0,T], \]
which implies $\eta^{(1)}=\eta^{(2)}$.

Changing the variable $t$ in \eqnu{yCy0eq0prf1} to $s$,
and then integrating from $0$ to $t$,
and changing the order of integration in the first term on the right
hand side, we have
\[ \arrb{l}\dsp
\int_0^t \eta(w,s)\,ds
\\ \dsp \phantom{\int_0^t}
= -\int_0^t \eta(w,u)\, \biggl( \int_u^t \pderiv{}{s}
 \exp (\dsp -\int_u^s w(\tilde{g}(u,v),v)\,dv ) ds\biggr)\, du
\\ \dsp \phantom{\int_0^t=}
 +\int_0^1 \pderiv{\rho}{z}(w,z)\, \biggl( \int_0^t \pderiv{}{s}
 \exp (\dsp -\int_0^s w(f(z,v),v)\,dv ) ds\biggr)\, dz
\\ \dsp \phantom{\int_0^t}
=
\int_0^t \eta(w,u)\,\biggl(1-
\exp (\dsp -\int_u^t w(\tilde{g}(u,v),v)\,dv )\biggr)\, du
\\ \dsp \phantom{\int_0^t=}
 -\int_0^1 \pderiv{\rho}{z}(w,z)\, \biggl(1-
 \exp (\dsp -\int_0^t w(f(z,v),v)\,dv )\biggr)\, dz,
\arre \]
which, with $\rho(w,0)=1$ and $\rho(w,1)=0$, proves \eqnu{yCy0eq0prf2}.

Combining \eqnu{yCy0eq0prf1} and \eqnu{yCy0eq0prf3},
together with 
$\dsp \pderiv{\rho}{z}(w,z)\le0$, 
$\rho(w,0)=1$ and $\rho(w,1)=0$,
we see that
\[
\int_W \eta(w,t)\,U_0(dw,0)\le 
R_w(T)\int_0^t \int_W \eta(w,u)\,U_0(dw,0)\,du +R_w(T).
\]
\cite[\S 3.8, Lemma 3.4]{Bressan} again implies that 
there exists $C(T)>0$, independent of $\tilde{g}$, 
such that $\dsp\int_W \eta(w,t)\,U_0(dw,0)\le C(T)$, $0\le t\le T$.
\QED\prfe
\corb
\cora{yCy0eq0prf4}
For $i=1,2$, let $\eta_{i}$ be $\eta$ in \lemu{yCy0eq0prf1}
with $g_i$ in place of $\tilde{g}$, respectively.
Then, for each $T>0$ there exists a positive constant $C(T)$ such that
\eqnb
\eqna{yCy0eq0prf4}
\int_W |\eta_{1}(w,t)-\eta_{2}(w,t)|\, U_0(dw,0)
\le C(T) \int_0^t \sup_{v\in [u,T]} |g_1(u,v)-g_2(u,v)|\,du.
\eqne
\DDD\core
\prfb
Put 
\[ \Delta \eta(t)=\int_W |\eta_{1}(w,t)-\eta_{2}(w,t)|\, U_0(dw,0) \]
and
\[ \Delta g(u)=\sup_{v\in [u,T]} |g_1(u,v)-g_2(u,v)|. \]
\lemu{yCy0eq0prf1}, in particular,
\eqnu{yCy0eq0prf1}, \eqnu{yCy0eq0prf2b}, and \eqnu{yCy0eq0prf3},
implies that
\[ \Delta \eta(t)
\le C_1(T)\int_0^t \Delta \eta(u)\, du
+C_2(T)\int_0^t \Delta g(u)\,du,\ t\in[0,T], \]
for each $T$ and for positive constants $C_i(T)$, $i=1,2$.
Hence
\[
\Delta \eta(t)\le C_2(T)\int_0^t e^{C_1(T)(t-s)}\Delta g(s)\,ds
\le C_2(T)e^{TC_1(T)}\int_0^t\Delta g(s)\,ds,
\]
which implies \eqnu{yCy0eq0prf4}.
\QED\prfe

\lemb
\lema{yCy0eq0}
There exists a unique $C^1$ function
$g:\ \{(s,t)\in \preals ^2\mid 0\le s\le t\}\to[0,1]$ such that
\eqnb
\eqna{yCy0eq0}
\arrb{l}\dsp
 g(s,t)=1+\int_W
 \int_{0}^1 \pderiv{\rho}{z}(w,z) \exp(-\int_0^t w(f(z,u),u)\,du)
 \,dz\,U_0(dw,0)
\\ \dsp \phantom{g(s,t)=1}
-\int_W \int_0^s \eta(w,u)
\exp( -\int_u^t w(g(u,v),v)\,dv )\,du\,U_0(dw,0),
\ \ 0\le s\le t.
\arre
\eqne
Here, $f(s,t)$ is defined in \eqnu{yCt0eq0}
and $\eta$ is the function given by \lemu{yCy0eq0prf1}
with $g$ in place of $\tilde{g}$.
\DDD\leme
\prfb
For $k\in\pintegers$, define 
a sequence of functions, $g_k$ and
$\eta_{k}$, inductively by
$\dsp g_{0}(s,t)=1$, $0\le s\le t$, and,
for $k\in\pintegers$,
$\dsp\eta_{k}$ the function $\eta$ in
\lemu{yCy0eq0prf1} with $g_k$ in place of $\tilde{g}$, and 
\eqnb
\eqna{yCy0eq0prf6}
\arrb{l}\dsp
 g_{k+1}(s,t)=1+\int_W
 \int_{0}^1 \pderiv{\rho}{z}(w,z) \exp(-\int_0^t (f(z,u),u)\,du) \,dz
 \, U_0(dw,0)
\\ \dsp \phantom{g_{k+1}(s,t)=1}
-\int_W \int_0^s \eta_{k}(u)\exp( -\int_u^t w(g_k(u,v),v)\,dv )\,du
 \, U_0(dw,0),
\ \ 0\le s\le t.
\arre
\eqne
Note that \eqnu{rhoasdistributionfunction} and $\eta_k(w,z)\ge0$
implies $g_k(s,t)\le 1$, and that
\eqnu{yCy0eq0prf2} and \eqnu{Burgersassump}, with $\eta_k(w,z)\ge0$ imply
\[ 1- g_k(s,t)\le \int_W \rho(w,0)\,U_0(dw,0)=1, \]
hence, 
$0\le g_k(s,t)\le 1$, implying that $\eta_k$ is well-defined.

Put 
$\Delta g_k=|g_{k+1}-g_k|$ and
$\Delta \eta_{k}=|\eta_{k+1}-\eta_{k}|$.
Repeating the arguments of \lemu{yCt0eq0} or \lemu{yCy0eq0prf1},
we see that \eqnu{yCy0eq0prf6} implies, with \eqnu{yCy0eq0prf2b},
\[ \Delta g_{k+1}(s,t)\le
\int_W \int_0^s \Delta \eta_{k}(w,u)\,du\,U_0(dw,0)
 +C_1(T)\int_0^s\biggl(\int_u^t \Delta g_k(u,v)\,dv\biggr) du,
\]
for $0\le s\le t\le T$, where $C_1(T)$ is a positive constant.
Putting $\dsp G_k(s)=\sup_{t\in[s,T]} \Delta g_{k}(s,t)$, we have, with
\coru{yCy0eq0prf4},
\[ \arrb{l}\dsp
G_{k+1}(s) \le C_2(T) \int_0^s \biggl(\int_0^u G_k(v)\,dv\biggr) du
 +T\,C_1(T)\int_0^s G_k(u)\,du
\\ \dsp \phantom{G_{k+1}(s)}
\le (C_2(T)+C_1(T))\,T\,\int_0^s G_k(u)\,du,
\arre \]
where $C_2(T)$ is a positive constant.
As in the proof of \lemu{yCt0eq0} or \lemu{yCy0eq0prf1},
this implies that the limit $g=\limf{k} g_k$ exists and is continuous.
Also, $0\le g_k(s,t)\le 1$ implies
\eqnb
\eqna{grange}
0\le g(s,t)\le 1,\ t\ge s\ge 0.
\eqne
Then 
$\eta=\limf{k} \eta_{k}$ also exist and are continuous,
and these functions satisfy
\eqnu{yCy0eq0prf1} with $g$ in place of $\tilde{g}$, and \eqnu{yCy0eq0}.
$C^1$ properties follow from the right hand side of \eqnu{yCy0eq0},
and uniqueness also follows as in the proof of \lemu{yCy0eq0prf1}.
\QED\prfe
\corb
\cora{yCy0eq0}
The following hold.
\eqnb
\eqna{yCt0eq0prf1}
f(y,0)=y,\ \ y\in[0,1].
\eqne
\eqnb
\eqna{fstrictonetoone}
\pderiv{f}{y}(y,t)>0,\ \pderiv{f}{t}(y,t)\ge 0,
\ \ (y,t)\in[0,1]\times\preals.
\eqne
\eqnb
\eqna{gprimegefprime}
\pderiv{g}{s}(s,t)\le 0,\ \pderiv{g}{t}(s,t)\ge 0,\ \ 0\le s\le t.
\eqne
\eqnb
\eqna{Ualongcharacteristicspre0s}
g(t,t)=0,\ \ t\ge0.
\eqne
\eqnb
\eqna{characteristicsconnection}
g(0,t)=f(0,t),\ t\ge0.
\eqne
\DDD\core
\prfb
The claims on $f$, \eqnu{yCt0eq0prf1} and \eqnu{fstrictonetoone},
are consequences of \eqnu{yCt0eq0} and the assumptions in\thmu{Burgers}.
The only perhaps less obvious claim is that the derivative of $f$ in $y$
cannot be $0$ in \eqnu{fstrictonetoone}, which follows from
\eqnu{rhoasdistributionfunction} and \eqnu{yCy0eq0prf3}, with
\[
\pderiv{f}{y}(y,t)
\ge -e^{-T\,R_w(T)} \int_W \pderiv{\rho}{y}(w,y)\,U_0(dw,0)
=  -e^{-T\,R_w(T)} \pderiv{}{y} U_0(\chrfcn{W},y)
=  e^{-T\,R_w(T)}>0.
\]

Differentiating \eqnu{yCy0eq0prf2} with $\tilde{g}$ replaced by $g$,
\eqnb
\eqna{yCy0eq0prf8}
\arrb{l} \dsp
\eta(w,t)-\int_0^t \eta(w,u)\,w(g(u,t),t)
\, \exp(-\int_u^t w(g(u,v),v)\,dv) \, du
\\ \dsp {}
= - \int_0^1 \pderiv{\rho}{z}(w,z) \,w(f(z,t),t)
 \,\exp(-\int_0^t w(f(z,s),s)\,ds)\,dz.
\arre
\eqne
Integrating \eqnu{yCy0eq0prf8} over $W$ with measure $U_0(dw,0)$,
and recalling that $\eta$ and $w\in W$ are non-negative, 
and using \eqnu{yCt0eq0} and \eqnu{fstrictonetoone},
\eqnb
\eqna{yCy0eq0prf9}
\arrb{l}\dsp
\int_W \eta(w,t)\, U_0(dw,0)\\ \dsp
\ge -\int_W \int_0^1 \pderiv{\rho}{z}(w,z) \,w(f(z,t),t)
 \,\exp(-\int_0^t w(f(z,s),s)\,ds)\,dz\, U_0(dw,0)\\ \dsp
=\pderiv{f}{t}(0,t)\ge 0.
\arre
\eqne
Differentiating \eqnu{yCy0eq0} by $s$, 
and using \eqnu{yCy0eq0prf9} we then have
\[
\pderiv{g}{s}(s,t)
=-\int_W \eta(w,s)\exp( -\int_s^t w(g(s,v),v)\,dv )\, U_0(dw,0)
\le 0. \]
Similarly, differentiating $g(s,t)$ by $t$ and using 
\eqnu{yCy0eq0} and \eqnu{fstrictonetoone},
\[
\pderiv{g}{t}(s,t)\ge \pderiv{f}{t}(0,t)\ge 0.
\]

The rest of the claims are obtained easily. Indeed,
\eqnu{Ualongcharacteristicspre0s} follows
from \eqnu{Burgersassump}, \eqnu{yCy0eq0} and \eqnu{yCy0eq0prf2},
and
\eqnu{characteristicsconnection} from \eqnu{yCy0eq0} and \eqnu{yCt0eq0}.
\QED\prfe
We are ready to define the characteristic curves $y=y_C(\gamma,t)$
for \eqnu{Burgers}.
For $\gamma=(y_0,t_0)\in\Gamma$ and $t\ge t_0$, put
\eqnb
\eqna{yCpre0}
y_C(\gamma,t):=\left\{\arrb{ll} f(y_0,t) & \mbox{ if } 
\gamma\in\Gamma_i\,,\ \mbox{ i.e., }t_0=0,
\\ g(t_0,t) & \mbox{ if } 
\gamma\in\Gamma_b\,,\ \mbox{ i.e., }y_0=0. \arre \right.
\eqne
Note that \eqnu{characteristicsconnection} implies that 
\eqnu{yCpre0} is well-defined on $(y_0,t_0)=(0,0)\in \Gamma_i\cap \Gamma_b$.
\lemu{yCt0eq0} and \lemu{yCy0eq0} imply continuity of $y_C(\gamma,t)$,
and $C^1$ property in $t$. (In fact, it is also $C^1$ in $(\gamma,t)$
except on $y=y_C((0,0),t)$.)
Also, \eqnu{yCt0eq0prf1} and \eqnu{Ualongcharacteristicspre0s}
imply the first equality in \eqnu{U0}.

Note also that, for each $t\ge0$, $y_C(\cdot,t):\ \Gamma_t\to [0,1]$
is surjective. In fact, $f$ and $g$ are continuous,
\eqnu{yCt0eq0prf1} and \eqnu{fstrictonetoone} imply
$f(1,t)\ge f(1,0)=1$.
These and \eqnu{Ualongcharacteristicspre0s} and 
\eqnu{characteristicsconnection} imply that $y_C$ is surjective:
\[
\{ y_C(\gamma,t) \mid \gamma\in\Gamma_t \}=[0,1].
\]

Note that \eqnu{fstrictonetoone} implies that
there exists a unique $C^1$, increasing, one-to-one onto inverse function 
$\hat{f}:\ [f(0,t),1]\to [0,1]$ of $f(y,t)$ with respect to $y$.
For $y< y_C(0,0,t)=f(0,t)=g(0,t)$ we define
$\hat{g}:\ [0,g(0,t)]\to [0,t]$ by
\eqnb
\eqna{hatg}
\hat{g}(y,t)=\inf\{ s\ge0 \mid g(s,t)=y \}.
\eqne
Since, as noted above, $g(\cdot,t):\ [0,t]\to [0,g(0,t)]$ is surjective,
$\hat{g}$ is well-defined, and 
Since $g$ is continuous, $g(\hat{g}(y,t),t)=y$.
Also \eqnu{gprimegefprime} implies that $\hat{g}(y,t)$ is 
non-increasing with respect to $y$.
Put
\eqnb
\eqna{hatgamma}
\hat{\gamma}(y,t)=\left\{\arrb{ll}
\dsp (\hat{f}(y,t),0)\in\Gamma_i & \mbox{ if } f(0,t)\le y\le 1, \\ 
\dsp (0,\hat{g}(y,t))\in\Gamma_b\cap\Gamma_t
 & \mbox{ if } 0\le y\le g(0,t). \arre\right.
\eqne
The definition implies
\eqnb
\eqna{hatgammayC}
y_C(\hat{\gamma}(y,t),t)=y,\ \ y\in[0,1],
\ \mbox{ and }\ 
\hat{\gamma}(y_C(\gamma,t),t)=\gamma,\ \ \gamma\in\Gamma_i.
\eqne
Note that the second equality may fail on $\gamma\in\Gamma_b$.

For $t\ge 0$, define a measure valued function
\[ \varphi(dw,\cdot,t):\ \Gamma_t\to\preals \]
as follows:
If $\gamma=(y_0,0)\in\Gamma_i$,
\eqnb
\eqna{Ualongcharacteristicsexplicit1}
\varphi(dw,\gamma,t):=
 - \int_{y_0}^1 \pderiv{\rho}{z}(w,z) 
\exp(-\int_0^t w(f(z,s),s)\,ds) \,dz\, U_0(dw,0),
\eqne
where $f$ is as in \lemu{yCt0eq0}, and 
if $\gamma=(0,t_0)\in\Gamma_b\cap\Gamma_t$,
\eqnb
\eqna{Ualongcharacteristicsexplicit2}
\ \ \ \arrb{l} \dsp
\varphi(dw,\gamma,t)
\\ \dsp {} :=
- \int_{0}^1 \pderiv{\rho}{z}(w,z) 
\exp(-\int_0^t w(f(z,s),s)\,ds) \,dz \,U_0(dw,0)
\\ \dsp \phantom{:=}
+ \int_0^{t_0} 
\eta(w,u)\exp( -\int_u^t w(g(u,v),v)dv )\,du\, U_0(dw,0),
\arre
\eqne
where, $f$ is as in \lemu{yCt0eq0}, and $\eta$ and $g$
are as in \lemu{yCy0eq0}.
Let
\[
\varphi(h,\gamma,t) :=\int _W h(s)\, \varphi(dw,\gamma,t)
\]
for a continuous bounded function $h$, 
$\gamma \in \Gamma$ and $t \in [0,\infty )$.
\prpb
\prpa{yCvarphi}
The following hold.
\eqnb
\eqna{yCt0eq02}
\eqna{yCpre}
y_C(\gamma,t)=1-\varphi(\chrfcn{W},\gamma,t)
:=1-\int_W \varphi(dw,\gamma,t),
\ \ \gamma\in\Gamma_t,\ t\ge 0.
\eqne
\eqnb
\eqna{Ualongcharacteristicspre0}
\varphi(dw,\gamma,t_0)=U_0(dw,y_0),
\ \ \gamma=(y_0,t_0)\in\Gamma.
\eqne
For bounded continuous $h:\ W\to\reals$ and $t>0$,
\eqnb
\eqna{densityalongcharacteristicsteq0}
\pderiv{\varphi}{t}(h,(y_0,0),t)=\int_W \int_{y_0}^1
 w(y_C((z,0),t),t)\pderiv{\varphi}{z}(h,(z,0),t)\,dz\,U_0(dw,0),
\ \ 0\le y_0\le 1,
\eqne
and
\eqnb
\eqna{densityalongcharacteristicsyeq0}
\arrb{l}\dsp
\pderiv{\varphi}{t}(h,(0,t_0),t)
=\pderiv{\varphi}{t}(h,(0,0),t)
-\int_W \int_0^{t_0} w(y_C((0,u),t),t)
\pderiv{\varphi}{u}(h,(0,u),t)\,du\,U_0(dw,0),
\\ \dsp
\ \ 0\le t_0\le t.
\arre
\eqne
\DDD\prpe
\prfb
The definitions \eqnu{yCpre0}, \eqnu{Ualongcharacteristicsexplicit1}
and \eqnu{Ualongcharacteristicsexplicit2}, with \lemu{yCt0eq0} and 
\lemu{yCy0eq0} imply \eqnu{yCpre}, and
\eqnu{Ualongcharacteristicspre0} follows from \eqnu{yCy0eq0prf2}, 
\eqnu{Ualongcharacteristicsexplicit1} and
\eqnu{Ualongcharacteristicsexplicit2}.
The definitions \eqnu{yCpre0} and \eqnu{Ualongcharacteristicsexplicit1}
imply that both hand sides of
\eqnu{densityalongcharacteristicsteq0} are equal to
\[ 
\int_W h(w)\int_{y_0}^1 \pderiv{\rho}{z}(w,z)\, w(f(z,t),t)\,
\exp(-\int_0^t w(f(z,s),s)\,ds)\,dz\,U_0(dw,0).
\]
Similarly, \eqnu{yCpre0} and \eqnu{Ualongcharacteristicsexplicit2}
imply that both hand sides of 
\eqnu{densityalongcharacteristicsyeq0} are equal to
\[ 
\pderiv{\varphi}{t}(h,(0,0),t)
- \int_W h(w) \int_0^{t_0} \eta(w,u)\,w(g(u,t),t)\,
\exp(-\int_u^t w(g(u,v),v)\,dv)\,du\,U_0(dw,0).
\]
\QED\prfe

For $(y,t)\in[0,1]\times\preals$ put
\eqnb
\eqna{U}
U(dw,y,t):=\varphi(dw,\hat{\gamma}(y,t),t)
=
\left\{\arrb{ll}\dsp
\varphi(dw,(\hat{f}(y,t),0),t) & f(0,t)\le y\le 1,
\\ \dsp
\varphi(dw,(0,\hat{g}(y,t)),t) & 0\le y\le g(0,t),
\arre\right.\eqne
where $\hat{\gamma}$ is defined in \eqnu{hatgamma}.
\thmb
\thma{yC}
It holds that
\eqnb
\eqna{Ualongcharacteristics}
\varphi(dw,\gamma,t)=U(dw,y_C(\gamma,t),t),
\ \ \gamma\in\Gamma_t,\ t\ge 0.
\eqne

Furthermore,
for bounded continuous function $h:\ W\to\reals$
$U(h,\cdot,\cdot):\ [0,1]\times\preals\to \preals$ is 
Lipschitz continuous in $(y,t)\in [0,1]\times[0,T]$ for any $T>0$,
and satisfies the second equality in \eqnu{U0}, 
\eqnu{phic}, \eqnu{yc}, and \eqnu{soliditycond}.
\DDD\thme
\prfb
For $\gamma\in\Gamma_i$, \eqnu{Ualongcharacteristics} follows from
\eqnu{hatgammayC}. The point is the case $\gamma\in \Gamma_b$,
where $y_C((0,s),t)=g(s,t)$, as a function of $s$, may fail to be one-to-one.
Suppose $g(s,t)=g(s',t)$ for some $s$ and $s'$ satisfying
$0\le s<s'\le t$.
Then \eqnu{yCy0eq0} and non-negativity of $\eta(w,u)$ implies
\[
\int_s^{s'} \eta(w,u)\exp(-\int_u^t w(g(u,v),v)dv)\,du=0,
\ U_0(dw,0) \mbox{--almost surely.}
\]
Hence \eqnu{Ualongcharacteristicsexplicit2} implies
$\varphi(dw,(0,s'),t)=\varphi(dw,(0,s),t)$.
On the other hand, the first equality of \eqnu{hatgammayC} implies
\[ y_C(\hat{\gamma}(y_C(\gamma,t),t),t)=y_C(\gamma,t),\ \gamma\in\Gamma_b. \]
Therefore,
$\varphi(dw,\hat{\gamma}(y_C(\gamma,t),t),t)=\varphi(dw,\gamma,t)$,
with which \eqnu{U} implies
\[
U(dw,y_C(\gamma,t),t)=\varphi(dw,\hat{\gamma}(y_C(\gamma,t),t),t)
=\varphi(dw,\gamma,t),
\]
so that \eqnu{Ualongcharacteristics} holds.

The Lipschitz continuity of $U(h,y,t)$ for
$f(0,t)\le y\le 1$, $0\le t\le T$ is obvious,
since the definitions \eqnu{U}, \eqnu{Ualongcharacteristicsexplicit1}, and
the definition of $\hat{f}$ stated just before \eqnu{hatg} imply
that $U(h,y,t)$ is $C^1$.
To prove the Lipschitz continuity of $U(h,y,t)$ for
$0\le g(0,t)=f(0,t)\le y\le 1$, $0\le t\le T$,
let $(y,t)$ and $(y',t')$ be $2$ points in this domain.
Use \eqnu{U} to decompose
\[
|U(h,y',t')-U(h,y,t)|
\le  |\varphi(h,\hat{\gamma}(y',t'),t')-\varphi(h,\hat{\gamma}(y',t'),t)|
+  |\varphi(h,\hat{\gamma}(y',t'),t)-\varphi(h,\hat{\gamma}(y,t),t)|.
\]
Since by definition \eqnu{Ualongcharacteristicsexplicit2} 
$\varphi(h,\gamma,t)$ is $C^1$ in $t$, the first term on the right hand
side is bounded by a global constant times $|t'-t|$.
To evaluate the second term,
let $M$ be such that $|h(w)|\le M$, $w\in W$, and
denote by $h_+$ and $h_-$ the positive and negative part of $h$,
respectively, so that $h=h_+-h_-$, $0\le h_{\pm}\le M$. 
Definitions \eqnu{hatgamma} and
\eqnu{Ualongcharacteristicsexplicit2},
and the non-negativity of $\eta$ imply
\[ \arrb{l}\dsp
|\varphi(h,\hat{\gamma}(y',t'),t)-\varphi(h,\hat{\gamma}(y,t),t)|
=
\biggl|
\int_{\hat{g}(y',t')}^{\hat{g}(y,t)} 
\int_W h(w)\,\eta(w,u)\exp( -\int_u^t w(g(u,v),v)dv )\, U_0(dw,0)\,du
\biggr|
\\ \dsp {}
\le
2M\,
\biggl|
\int_{\hat{g}(y',t')}^{\hat{g}(y,t)} 
\int_W \chrfcn{W}(w)\,\eta(w,u)\exp( -\int_u^t w(g(u,v),v)dv )
\, U_0(dw,0)\,du
\biggr|,
\arre \]
which, with with \eqnu{hatgamma}, \eqnu{Ualongcharacteristicsexplicit2},
and \eqnu{yCt0eq02}, is equal to
\[ 
2M\,
|y_C(\hat{\gamma}(y',t'),t) - y_C(\hat{\gamma}(y,t),t)|.
\]
This with \eqnu{hatgammayC} implies
\[ \arrb{l} \dsp
|\varphi(h,\hat{\gamma}(y',t'),t)-\varphi(h,\hat{\gamma}(y,t),t)|
\\ \dsp {}
\le 2M\,|y_C(\hat{\gamma}(y',t'),t') - y_C(\hat{\gamma}(y',t'),t)|
+ 2M\,|y_C(\hat{\gamma}(y',t'),t') - y_C(\hat{\gamma}(y,t),t)|
\\ \dsp \phantom{\le}
= 2M\,|y_C(\hat{\gamma}(y',t'),t') - y_C(\hat{\gamma}(y',t'),t)|
+ 2M\,|y'-y|.
\arre\]
Since $y_C(\gamma,t)$ is $C^1$ in $t$, we have the global Lipschitz
continuity.

The property \eqnu{soliditycond} follows from \eqnu{yCt0eq02} 
and \eqnu{hatgammayC}.
The second equality in \eqnu{U0} then follows from
\eqnu{Ualongcharacteristicspre0}, \eqnu{Ualongcharacteristics},
\eqnu{Ualongcharacteristics}, and the first equality in \eqnu{U0}.
(Note that the first equality in \eqnu{U0} and other claims in
\thmu{Burgers} for $y_C$ is proved below \eqnu{yCpre0}.)

To prove \eqnu{phic} for $(y_0,t_0)\in\Gamma_i$, namely, for $t_0=0$,
use \eqnu{prfofHDL3}, \eqnu{U}, \eqnu{Ualongcharacteristics},
and \eqnu{Ualongcharacteristicsexplicit1},
and change the order of integration, to find
\[ \arrb{l}\dsp
-V(h,y_C((y_0,0),t),t)
\\ \dsp {}
=-\int_W h(w)\, w(y_C((y_0,0),t),t)\,\varphi(dw,(y_0,0),t)
\\ \dsp \phantom{=}
+\int_W h(w)\, \int_{y_C((y_0,0),t)}^1 \pderiv{w}{z}(z,t)\,\biggl(
 \int_{\hat{f}(z,t)}^1 \pderiv{\rho}{z'}(w,z') 
\exp(-\int_0^t w(f(z',s),s)\,ds) \,dz'\biggr)\,dz\,U_0(dw,0)
\\ \dsp {}
=-\int_W h(w)\,w(y_C((y_0,0),t),t)\,\varphi(dw,(y_0,0),t)
\\ \dsp \phantom{=}
+ \int_W h(w)\,\int_{y_0}^1 
\biggl(\int_{y_C((y_0,0),t)}^{y_C((z',0),t)} 
\pderiv{w}{z}(z,t)dz\biggr)
\pderiv{\rho}{z'}(w,z') \exp(-\int_0^t w(f(z',s),s)\,ds)
\, dz'\,U_0(dw,0),
\arre \]
which, with the definition \eqnu{Ualongcharacteristicsexplicit1},
is equal to $\dsp \pderiv{\varphi}{t}(h,\gamma,t)$.
Integrating from $t_0$ to $t$ and using
\eqnu{Ualongcharacteristics} and \eqnu{Ualongcharacteristicspre0},
we have \eqnu{phic}.

To prove \eqnu{phic} for $(y_0,t_0)\in\Gamma_b$, namely, for $y_0=0$,
first decompose the integration range in \eqnu{prfofHDL3} with
$y=y_C((0,t_0),t)$ as
\[ [y_C((0,t_0),t),1]= [g(t_0,t),g(0,t)]\cup[f(0,t),1], \]
then use the definitions \eqnu{U} and 
\eqnu{Ualongcharacteristicsexplicit1} or
\eqnu{Ualongcharacteristicsexplicit2},
and change the order of integration, to find
\[ \arrb{l}\dsp
-V(h,y_C((0,t_0),t),t)
\\ \dsp {}
=-\int_W h(w)\,w(y_C((0,t_0),t),t)\varphi(dw,(0,t_0),t)
\\ \dsp \phantom{=}
-\int_W h(w)\,\int_{g(t_0,t)}^{g(0,t)} \pderiv{w}{z}(z,t)\,\biggl(
- \int_0^1 \pderiv{\rho}{z'}(w,z') 
\exp(-\int_0^t w(f(z',s),s)\,ds) \,dz'
\\ \dsp 
\phantom{=-\int_W h(w)\,\int_{g(t_0,t)}^{g(0,t)} 
\pderiv{w}{z}(z,t)\,\biggl(}
+\int_0^{\hat{g}(z,t)}
 \eta(w,u)\exp( -\int_u^t w(g(u,v),v)dv )\,du
\biggr)\,dz\,U_0(dw,0)
\\ \dsp \phantom{=}
+\int_W h(w)\,\int_{f(0,t)}^1 \pderiv{w}{z}(z,t)\,\biggl(
 \int_{\hat{f}(z,t)}^1 \pderiv{\rho}{z'}(w,z') 
\exp(-\int_0^t w(f(z',s),s)\,ds) \,dz'\biggr)\,dz\,U_0(dw,0)
\\ \dsp {}
=-\int_W h(w)\,w(y_C((0,t_0),t),t)\varphi(dw,(0,t_0),t)
\\ \dsp \phantom{=}
+\int_W h(w)\,(w(g(0,t),t)-w(g(t_0,t),t))\,
 \int_0^1 \pderiv{\rho}{z'}(w,z') 
\exp(-\int_0^t w(f(z',s),s)\,ds) \,dz'\,U_0(dw,0)
\\ \dsp \phantom{=}
-\int_W h(w)\,\int_0^{t_0}
 \eta(w,u)\exp( -\int_u^t w(g(u,v),v)dv )\,
\biggl(\int_{g(t_0,t)}^{g(u,t)} \pderiv{w}{z}(z,t)\,dz\biggr)
\,du\,U_0(dw,0)
\\ \dsp \phantom{=}
+ \int_W h(w)\,\int_{0}^1 \pderiv{\rho}{z'}(w,z') 
\exp(-\int_0^t w(f(z',s),s)\,ds) 
\,\biggl(\int_{f(0,t)}^{f(z',t)} \pderiv{w}{z}(z,t)\,dz\biggr)
\,dz'\,U_0(dw,0).
\arre \]
Using 
\eqnu{Ualongcharacteristicsexplicit2}, this further is simplified as
\[ \arrb{l}\dsp
\int_W h(w)\,\biggl(-\int_0^{t_0} \eta(u) w(g(u,t),t)
\exp( -\int_u^t w(g(u,v),v)dv )\,du
\\ \dsp \phantom{\int_W h(w)\,\biggl(}
+ \int_0^1 \pderiv{\rho}{z'}(w,z')  w(f(z',t),t)
\exp(-\int_0^t w(f(z',s),s)\,ds) \,dz'\biggr)\,U_0(dw,0),
\arre\]
which, by using \eqnu{Ualongcharacteristicsexplicit2}, 
is seen to be equal to $\dsp \pderiv{\varphi}{t}(h,\gamma,t)$.
Integrating from $t_0$ to $t$ and using
\eqnu{Ualongcharacteristics} and \eqnu{Ualongcharacteristicspre0},
we have \eqnu{phic}.

Substituting $h=\chrfcn{W}$ in \eqnu{phic},
and using \eqnu{soliditycond}, \eqnu{yCpre} and \eqnu{Ualongcharacteristics}, 
we have \eqnu{yc}.
\QED\prfe

To complete a proof of \thmu{Burgers},
it only remains to prove uniqueness.
Besides the pair $y_C$ and $U$ which we constructed and proved so far to 
satisfy the properties stated in \thmu{Burgers},
assume that there are another such pair $\tilde{y}_C$ and $\tilde{U}$.
For $T>0$, let $L(T)>0$ be such that
\[ \arrb{l}\dsp
\max\{ |U(h,y,t)-U(h,y',t')|,\ |\tilde{U}(h,y,t)-\tilde{U}(h,y',t')| \}
\\ \dsp{}
\le L(T)\, \norm{(y,t)-(y',t')},\ \ (y,t),(y',t')\in[0,1]\times[0,T],
\ h:\ W\to[-1,1];\mbox{ conti.}.
\arre \]
Put 
\[ I(t)=\sup_{h:\ W\to[-1,1]; \mbox{ conti.}}
\sup_{y\in [0,1]} |U(h,y,t)-\tilde{U}(h,y,t)| \]
and
\[
J(t)=
\sup_{\gamma\in\Gamma_t}|\tilde{y}_C(\gamma,t)-y_C(\gamma,t)|.
\]
Then \eqnu{U0} and its correspondence for $\tilde{U}$ imply
$I(0)=0$.
Since $y_C(\cdot,t):\ \Gamma_t\to [0,1]$ is onto,
\eqnb
\eqna{yCsurjI}
 I(t)=\sup_{h:\ W\to[-1,1]; \mbox{ conti.}} \sup_{\gamma\in \Gamma_t}
 |U(h,y_C(\gamma,t),t)-\tilde{U}(h,y_C(\gamma,t),t)|. 
\eqne
Note also that since $\tilde{U}(dw,y,t)$ is, by assumption,
a non-negative measure, for $h$ with $|h(w)|\le 1$, $w\in W$,
we have 
\[
\tilde{U}(h,y,t)\le \tilde{U}(\chrfcn{W},y,t)=1-y\le 1,
\]
where we also used \eqnu{soliditycond}.

It holds that
\[ I(t)\leq L(T)J(t)\]
Subtracting 
\eqnb
\eqna{phictilde}
\tilde{U}(h,\tilde{y}_C(\gamma,t),t)
=U_0(h,y_0)-\int_{t_0}^t \tilde{V}(h,\tilde{y}_C(\gamma,s),s)\,ds,
\eqne
from \eqnu{phic}, and using \eqnu{yCsurjI},
\eqnu{prfofHDL3} and \eqnu{yCy0eq0prf3},  we have
\[ \arrb{l}\dsp
|\tilde{y}_C(\gamma,t)-y_C(\gamma,t)|
=
|\tilde{U}(\chrfcn{W},\tilde{y}_C(\gamma,t),t)
-U(\chrfcn{W},y(\gamma,t),t)|
\\ \dsp {}
\le
2R_w(T)\,\int_{t_0}^t J(s)\,ds
+R_w(T)\,\int_{t_0}^t I(s)\,ds
+R_w(T)L(T)\,\int_{t_0}^t J(s)\,ds.
\arre \]
Therefore,
\[ 
J(t)\le 2R_w(T)\,\int_{t_0}^t J(s)\,ds+R_w(T)\,\int_{t_0}^t I(s)ds
+R_w(T)L(T)\,\int_{t_0}^t J(s)\,ds.
\]
Then,
\[ 
I(t)
\le L(T)\,
\biggl(2R_w(T)\,\int_{t_0}^t J(s)\,ds+R_w(T)\,\int_{t_0}^t I(s)ds
+R_w(T)L(T)\,\int_{t_0}^t J(s)\,ds \biggr),
\]
so that if we put $K(t)=\max\{ I(t), J(t) \}$ then
there exists $C(T)$ such that
\[ K(t)\le C(T)\int_{t_0}^t K(s)\,ds, \ K(t_0)=0, \]
implying $K(t)=0$.
Hence $\tilde{U}=U$ and $\tilde{y}_C=y_C$.

This completes a proof of \thmu{Burgers}.

\section{Proof of \protect\thmu{HDLposdep}.}
\seca{profofHDL}

Let $\Gamma$ be as in \eqnu{spacetimeboundary}.
To simplify the notation, 
for $\gamma=(y_0,t_0)\in \Gamma$,
we will write $y_C((y_0,t_0),t)$ defined in \eqnu{yCpre0}
as $y_C(y_0,t_0,t)$.

We first prepare a random variable which converges as $N\to\infty$ 
to $y_C(t)=y_C(y_0,t_0,t)$, for $(y_0,t_0)\in \Gamma$, $t\geq t_0$;
\eqnb
\eqna{yCN}
Y\ssN_C(y_0,t_0,t)=y_0+\frac1N\sum_{i;\ X\ssN_i(t_0) \ge Ny_0+1}
 \chrfcn{J\ssN_i(t_0,t)}\,,\ \ (y_0,t_0)\in[0,1)\times [0,\infty ),\ t>t_0\,,
\eqne
where $J\ssN_i$ is defined in \eqnu{Jtopist}.
In particular, if we put, as an analogue to \eqnu{YN},
\eqnb
\eqna{yNinit}
y\ssN_i=\frac1N(x\ssN_i-1),\ \ i=1,2,\ldots,N,
\eqne
then \eqnu{Jtopist} and \eqnu{Jstepdown} imply
\eqnb
\eqna{JYcboundary}
Y\ssN_i(t)\ge Y\ssN_C(y_0,t_0,t)\ \ \Leftrightarrow
\ \ Y\ssN_i(t_0)\ge y_0\ \mbox{ and }\ J\ssN_i(t_0,t)\ \mbox{ does not hold}.
\eqne
Hence, we have
\eqnb
\eqna{YNC2}
\arrb{l}
Y\ssN_C(y_0,t_0,t) \\[2mm] \dsp
= y_0 + \frac1N \sum_{i} \int _{s\in (t_0,t]} \int _{\xi \in [0,\infty )} \chrfcn{Y\ssN_i(s-)\geq Y\ssN_C(y_0,t_0,s-)} \chrfcn{\xi \in [0,w_i(Y\ssN_i(s-),s))}\nu _i (d\xi ds).
\arre
\eqne
For the spatially homogeneous case, 
$Y\ssN_A(t_0,t)$ in \cite{timedep} is equal to $Y\ssN_C(0,t-t_0,t)$,
$Y\ssN_B(y_0,t)$ to $Y\ssN_C(y_0,0,t)$,
and $Y\ssN_C(t)$ in \cite{timedep} is equal to $Y\ssN_C(0,0,t)$ 
of \eqnu{yCN}.

Let $\Gamma$ be as in \eqnu{spacetimeboundary}.
Let $(y_0,t_0)\in \Gamma$, $t\geq t_0$. 
The definition \eqnu{distributionfcnexplicit} and the properties
\eqnu{Jtopist}, \eqnu{Jstepdown}, and \eqnu{JYcboundary}
imply that for $B\in {\mathscr B}(W)$, $U\ssN(B,Y\ssN_C(y_0,t_0,t),t)$ as a function of $t$
changes its value if and only if $J\ssN_i(t_0,t)$ occurs for some
$i$ satisfying $y\ssN_i\ge y_0$ and $w_i\in B$.
Therefore,  for $B\in {\mathscr B}(W)$
\[ \arrb{l}\dsp
U^{(N)}(B,Y^{(N)}_C(y_0,t_0,t),t) - U^{(N)}(B,y_0,t_0) \\[2mm] \dsp
= -\frac 1N \sum_{i;\ w_i \in B} \int _{s\in (t_0,t]} \int _{\xi \in [0,\infty)} \chrfcn{Y\ssN_i(s-)\geq Y\ssN_C(y_0,t_0,s-)} \chrfcn{\xi \in [0,w_i(Y_i^{(N)}(s-),s))} \nu _i(d\xi ds).
\arre \]

In analogy to \eqnu{prfofHDL3} define for $B\in {\mathscr B}(W)$
\eqnb
\eqna{prfofHDL2}
V\ssN (B,y,t)=\int _B w(y,t)\,U\ssN (dw,y,t)
+\int_y^1 \int _B \pderiv{w}{z}(z,t)\,U\ssN (dw,z,t)\,dz.
\eqne
By definition \eqnu{distributionfcnexplicit}, for $B\in {\mathscr B}(W)$
\eqnb
\eqna{prfofHDL4}
\eqna{prfofHDL5}
\arrb{l}\dsp
V\ssN (B,y,t)
= \frac1N \sum_{m\ge Ny+1} \sum_{i;\ w_i \in B} w_i(\frac{m-1}N,t)\,
\chrfcn{X\ssN_i(t-)=m}\,\\ \dsp
\phantom{V\ssN (B,y,t)} =\frac1N \sum_{i;\ w_i \in B,\, Y^{(N)}_i(t-)\ge y} w_i(Y_i^{(N)}(t-),t).
\arre
\eqne

Denote the compensated Poisson process by
\eqnb
\eqna{compensatedPp}
\tilde{\nu}_i(d\xi ds)=\nu _i(d\xi ds)-d\xi ds,
\eqne
and put for $B\in {\mathscr B}(W)$
\eqnb
\eqna{MU}
\arrb{l}
M^{(N)}_{U}(B,y_0,t_0,t) \\[2mm]
\dsp = -\frac 1N  \sum_{i;\ w_i \in B}
 \int _{s\in (t_0,t]} \int _{\xi \in [0,\infty)}  \chrfcn{Y\ssN_i(s-)\geq Y\ssN_C(y_0,t_0,s-)} \chrfcn{\xi \in [0,w_i(Y_i^{(N)}(s-),s))} \tilde{\nu}_i(d\xi ds).
\arre
\eqne
Then, we have for $B\in {\mathscr B}(W)$
\[
\arrb{l}
\dsp U^{(N)}(B,Y^{(N)}_C(y_0,t_0,t),t) \\
\dsp = U^{(N)}(B,y_0,t_0) + M^{(N)}_{U}(B,y_0,t_0,t) \\[2mm] \dsp
\quad - \frac 1N  \sum_{i;\ w_i \in B} \int _{s\in (t_0,t]} \int _{\xi \in [0,\infty)}  \chrfcn{Y\ssN_i(s-)\geq Y\ssN_C(y_0,t_0,s-)} \chrfcn{\xi \in [0,w_i(Y_i^{(N)}(s-),s))}d\xi ds \\
\dsp = U^{(N)}(B,y_0,t_0) + M^{(N)}_{U}(B,y_0,t_0,t)  \\[2mm] \dsp
\quad - \frac 1N  \sum_{i;\ w_i \in B} \int _{t_0}^t w_i \left( Y_i^{(N)}(s-),s\right)  \chrfcn{Y\ssN_i(s-)\geq Y\ssN_C(y_0,t_0,s-)} ds \\
\dsp = U^{(N)}(B,y_0,t_0) + M^{(N)}_{U}(B,y_0,t_0,t) - \int _0^t V^{(N)}(B,Y^{(N)}_C(y_0,t_0,s),s) ds.
\arre
\]
Combining this equality with \eqnu{phic},
we have for $B\in {\mathscr B}(W)$
\eqnb
\eqna{prfofHDL7U}
\arrb{l}\dsp
U\ssN (B,Y\ssN_C(y_0,t_0,t),t)- U(B,y_C(y_0,t_0,t),t)
\\[2mm] \dsp =U\ssN (B,y_0,t_0)- U(B,y_0,t_0) + M\ssN_{U}(B,y_0,t_0,t)
\\ \dsp
\quad -\int_{t_0}^t 
\biggl(V\ssN (B,Y\ssN_C(y_0,t_0,s),s) -V(B,y_C(y_0,t_0,s),s)\biggr)
\,ds.
\arre
\eqne
Put
\eqnb
\eqna{prfofHDL9}
\arrb{l}\dsp
W\ssN (t)=
\sup _{(y_0, t_0)\in \Gamma ; t_0\leq t}\ \sup_{s\in [t_0,t]}|Y\ssN_C(y_0,t_0,s)-y_C(y_0,t_0,s)|
\\ \dsp \phantom{W\ssN(t)=}
\;\vee\; \sup _{(y_0, t_0)\in \Gamma ; t_0\leq t}\ \sup_{s\in [t_0,t]}
||U\ssN (\cdot ,Y\ssN_C(y_0,t_0,s),s) \\ \dsp
\phantom{W\ssN _h(t)= \;\vee\;  \sup _{(y_0, t_0)\in \Gamma ; t_0\leq t}\ \sup_{s\in [t_0,t]} -\int _W h(w)}-U(\cdot ,y_C(y_0,t_0,s),s)||_{\rm{var}}
\\ \dsp \phantom{W\ssN(t)=}
\;\vee\;
\sup _{B\in {\mathscr B}(W)} \sup _{(y_0, t_0)\in \Gamma ; t_0\leq t}\ \sup_{s\in [t_0,t]}
|V\ssN (B,Y\ssN_C(y_0,t_0,s),s) \\ \dsp
\phantom{W\ssN _h(t)= \;\vee\; \sup _{(y_0, t_0)\in \Gamma ; t_0\leq t}\ \sup_{s\in [t_0,t]} -\int _W h(w)}-V(B,y_C(y_0,t_0,s),s)|.
\arre
\eqne
Since for all $z\in[0,1)$ there exist $(y_0,t_0)\in \Gamma$ such that $y_C(y_0,t_0,t)=z$, for $B\in {\mathscr B}(W)$
\[
\arrb{l}\dsp
\sup _{z\in [0,1]} |U\ssN (B,z,t)-U(B,z,t)|\\
\dsp
\leq \sup _{(y_0, t_0)\in \Gamma ; t_0\leq t}\ \sup_{s\in [t_0,t]} |U\ssN (B,y_C(y_0,t_0,s),s)-U(B,y_C(y_0,t_0,s),s)|\\
\dsp
\leq \sup _{(y_0, t_0)\in \Gamma ; t_0\leq t}\ \sup_{s\in [t_0,t]} |U\ssN (B,Y\ssN_C(y_0,t_0,s),s) -U(B,y_C(y_0,t_0,s),s)| \\
\dsp \quad + \sup _{(y_0, t_0)\in \Gamma ; t_0\leq t}\ \sup_{s\in [t_0,t]} |U\ssN (B,Y\ssN_C(y_0,t_0,s),s) -U\ssN (B,y_C(y_0,t_0,s),s)|.
\arre
\]
Hence,
\eqnb
\eqna{prfofHDL1001}
\arrb{l}\dsp
\sup _{z\in [0,1]} |U\ssN (B,z,t)-U(B,z,t)|\\
\dsp
\leq W\ssN (t) + \sup _{(y_0, t_0)\in \Gamma ; t_0\leq t}\ \sup_{s\in [t_0,t]} |U\ssN (B,Y\ssN_C(y_0,t_0,s),s) \\ \dsp
\phantom{\leq W\ssN (t) + \sup _{(y_0, t_0)\in \Gamma ; t_0\leq t}\ \sup_{s\in [t_0,t]} }-U\ssN (B,y_C(y_0,t_0,s),s)|.
\arre
\eqne
By \eqnu{distributionfcnexplicit} it holds that
\eqnb
\eqna{prfofHDL1002}
\arrb{l}\dsp
\sup _{(y_0, t_0)\in \Gamma ; t_0\leq t}\ \sup_{s\in [t_0,t]}|U\ssN (B,Y\ssN_C(y_0,t_0,s),s) -U\ssN (B,y_C(y_0,t_0,s),s)| \\
\dsp \hspace{4cm}\leq \sup _{(y_0, t_0)\in \Gamma ; t_0\leq t}\ \sup_{s\in [t_0,t]}|Y\ssN_C(y_0,t_0,s)-y_C(y_0,t_0,s)| +\frac{1}{N}.
\arre
\eqne
By \eqnu{prfofHDL1001} and \eqnu{prfofHDL1002} we obtain
\eqnb
\arrb{l}\dsp
\sup _{z\in [0,1]} |U\ssN (B,z,t)-U(B,z,t)| \leq 2W\ssN (t) +\frac{1}{N}.
\arre
\eqne
This implies
\eqnb
\eqna{prfofHDL9-2}
\arrb{l}\dsp
\sup _{z\in [0,1]} ||U\ssN (\cdot ,z,t)-U(\cdot ,z,t)||_{\rm{var}} \leq 4W\ssN (t) +\frac{2}{N}.
\arre
\eqne
By \eqnu{prfofHDL7U} and \eqnu{prfofHDL9}, we have for $B\in {\mathscr B}(W)$
\eqnb
\eqna{prfofHDL7U2}
\arrb{l}\dsp
\left| U\ssN (B,Y\ssN_C(y_0,t_0,t),t) - U(B,y_C(y_0,t_0,t),t) \right|
\\[2mm] \dsp \leq \left| U\ssN(B,y_0,t_0) - U(B,y_0,t_0)\right| + \left| M\ssN_{U}(B,y_0,t_0,t)\right| + \int _{t_0}^t W\ssN(s) ds.
\arre
\eqne
Similarly, combining \eqnu{prfofHDL2} with \eqnu{prfofHDL3}, we have
\[
\arrb{l}\dsp
V\ssN (B,Y\ssN_C(y_0,t_0,t),t) -V(B,y_C(y_0,t_0,t),t)
\\[2mm] \dsp {}
=
\int _B w(Y\ssN_C(y_0,t_0,t),t)\,
\bigl(U\ssN (dw,Y\ssN_C(y_0,t_0,t),t) -U(dw,y_C(y_0,t_0,t),t)\bigr)
\\[2mm] \dsp \phantom{=}
+\int _B \bigl(w(Y\ssN_C(y_0,t_0,t),t)-w(y_C(y_0,t_0,t),t)\bigr)
\,\,U(dw,y_C(y_0,t_0,t),t)
\\[2mm] \dsp \phantom{=}
+\int_{Y\ssN_C(y_0,t_0,t)}^1 \int _B \pderiv{w}{z}(z,t)\,
\bigl(U\ssN (dw,z,t)-U(dw,z,t)\bigr)\,dz
\\[2mm] \dsp \phantom{=}
-\int_{y_C(y_0,t_0,t)}^{Y\ssN_C(y_0,t_0,t)} \int _B \pderiv{w}{z}(z,t)\,
U(dw,z,t)\,dz.
\arre
\]
Hence, using this estimate, \eqnu{rw}, \eqnu{prfofHDL9-2} 
and the fact that $0\leq U\ssN \leq 1$, we have for $B\in {\mathscr B}(W)$
\eqnb
\eqna{prfofHDL7V2}
\arrb{l}\dsp
\left| V\ssN (B,Y\ssN_C(y_0,t_0,t),t)
-V(B,y_C(y_0,t_0,t),t)\right| \\[2mm]
\leq 
 R_w(T) || U\ssN (\cdot ,Y\ssN_C(y_0,t_0,t),t)
-U(\cdot ,y_C(y_0,t_0,t),t)||_{\rm{var}}
\\[2mm] \dsp \quad + 2 R_w(T) \left| Y\ssN_C(y_0,t_0,t)-y_C(y_0,t_0,t)\right|
+ R_w(T)\left( 4\int _{t_0+}^t W\ssN(s) ds +\frac{2}{N}\right) .
\arre
\eqne
To estimate $Y\ssN_C(y_0,t_0,t)-y_C(y_0,t_0,t)$, by using \eqnu{YNC2}, \eqnu{prfofHDL4} and \eqnu{MU} calculate
\eqnb
\eqna{YNC}
Y\ssN_C(y_0,t_0,t)=y_0+ M^{(N)}_{U}(W,y_0,t_0,t)
+ \int_{t_0}^t V\ssN(W,Y\ssN_C(y_0,t_0,s),s)\,ds.
\eqne
Combining with \eqnu{yc},
\eqnb
\eqna{prfofHDL7Y}
\arrb{l}
Y\ssN_C(y_0,t_0,t)-y_C(y_0,t_0,t)  \\[2mm] \dsp
= M^{(N)}_{U}(W,y_0,t_0,t) + \int_{t_0}^t
\left[ V\ssN( W,Y\ssN_C(y_0,t_0,s),s)-V(W,y_C(y_0,t_0,s),s)\right]
\,ds.
\arre
\eqne
Hence, by \eqnu{prfofHDL9} we have
\eqnb
\eqna{prfofHDL7Y2}
\left| Y\ssN_C(y_0,t_0,t)-y_C(y_0,t_0,t) \right|
\dsp \leq |M^{(N)}_{U}(W,y_0,t_0,t)| + \int _{t_0}^t W\ssN(s) ds.
\eqne
By \eqnu{rw}, \eqnu{prfofHDL7U2}, \eqnu{prfofHDL7V2}, and \eqnu{prfofHDL7Y2},we obtain for $B\in {\mathscr B}(W)$
\eqnb
\eqna{prfofHDL7V3}
\arrb{l}\dsp
\left| V\ssN(B,Y\ssN_C(y_0,t_0,t),t)-V(B,y_C(y_0,t_0,t),t)\right| \\[2mm]
\dsp \leq 
 R_w(T) || U\ssN(\cdot ,y_0,t_0)-U(\cdot ,y_0,t_0)||_{\rm{var}} + 7 R_w(T) \int _{t_0}^t W\ssN(s) ds \\[2mm] \dsp \quad
 + 3R_w(T)\left| M\ssN_{U}(B,y_0,t_0,t)\right| +  \frac{2R_w(T)}{N} .
\arre
\eqne
Because of \eqnu{prfofHDL9}, \eqnu{prfofHDL7U2}, \eqnu{prfofHDL7Y2}, and \eqnu{prfofHDL7V3}, we have
\[
\arrb{l}
\dsp W^{(N)}(t)
\leq C_1 \sup _{(y_0, t_0)\in \Gamma}|| U\ssN(\cdot,y_0,t_0)-U(\cdot,y_0,t_0)||_{\rm{var}}  \\[2mm]
\phantom{W^{(N)}(t)
\leq}\dsp + \left( 1+R_w(T)\right) \sup_{B\in {\mathscr B}(W)} \sup _{(y_0, t_0)\in \Gamma ; t_0\leq t}\left| M\ssN_{U}(B,y_0,t_0,t)\right| + C_2\ \int _0^t W\ssN(s) ds + \frac{2R_w(T)}{N}
\arre
\]
where $C_1$ and $C_2$ are constants depending on $A, T, R_w(T)$.
Hence, Gronwall's inequality implies
\eqnb
\eqna{prfofHDL7W}
\arrb{l}
\dsp \sup _{t\in [0,T]}W^{(N)}(t) \\[2mm]
\dsp
\leq e^{C_2T} \left[ C_1 \sup _{(y_0, t_0)\in \Gamma}|| U\ssN(\cdot,y_0,t_0)-U(\cdot,y_0,t_0)||_{\rm{var}} \right. \\[2mm]
\dsp \hspace{2cm} \left. + \left( 1+R_w(T)\right) \sup_{B\in {\mathscr B}(W)} \sup _{(y_0, t_0)\in \Gamma ; t_0\leq t}\left| M\ssN_{U}(B,y_0,t_0,t)\right|  + \frac{2R_w(T)}{N} \right].
\arre
\eqne

By the definition of $\Gamma$, we have
\eqnb
\eqna{prfofHDL1006}
\arrb{l} \dsp
E\left[ \sup_{B\in {\mathscr B}(W)} \sup _{(y_0, t_0)\in \Gamma , t_0\leq T}\ \sup _{t\in [t_0,T]}\left| M\ssN_{U}(B,y_0,t_0,t)\right| ^2\right] \\[2mm] \dsp
\leq E\left[ \sup_{B\in {\mathscr B}(W)} \sup _{y_0 \in [0,1)}\ \sup _{t\in [0,T]}\left| M\ssN_{U}(B,y_0,0,t)\right| ^2\right]  \\[2mm] \dsp
\quad + E\left[ \sup_{B\in {\mathscr B}(W)} \sup _{0\leq t_0 \leq t \leq T}\left| M\ssN_{U}(B,0,t_0,t)\right| ^2\right] .
\arre
\eqne
Note that $Y_i^{(N)}(t) \in \{ 0,1/N,\dots (N-1)/N\}$ for $t\in [0,T]$ and $i=1,2,\dots ,N$, $Y_C^{(N)}(y_0,t_0,\cdot )$ is a process of pure jumps by $1/N$ and for each $k=1,2,\dots ,N$, $Y_C^{(N)}(y_0,t_0,\cdot ) -y_0$ is independent of $y_0$ as long as $y_0\in ({k-1}/N,k/N]$. Also, note that $\{ M\ssN_{U}(B,\cdot ,0,\cdot ); B\in {\mathscr B}(W)\} = \{ M\ssN_{U}(B,\cdot ,0,\cdot ); B\in 2^{\{ w_i; i=1,2,\dots ,N\}}\}$.
Hence, by \eqnu{MU}
\[
\arrb{l} \dsp
E\left[ \sup_{B\in {\mathscr B}(W)}\sup _{y_0 \in [0,1)}\ \sup _{t\in [0,T]}\left| M\ssN_{U}(B,y_0,0,t)\right| ^2\right] \\[2mm] \dsp
= \frac 1{N^2} E\left[ \sup_{B\in {\mathscr B}(W)} \sup _{y_0 \in [0,1)}\ \sup _{t\in [0,T]}\left| \sum_{i;\, w_i\in B} \int _{s\in (0,t]}\int _{\xi \in [0,\infty )} \right. \right. \\[2mm] \dsp
\left.\left.\phantom{= \frac 1{N^2} \sum_{i;\, w_i\in B} \int _{s\in (0,t]}\int _{\xi \in [0,\infty )}}
\chrfcn{Y_i^{(N)}(s-)\geq Y_C^{(N)}(y_0,0,s-)} \chrfcn{\xi\in [0,w_{i}(Y_i^{(N)}(s-),s))} \tilde{\nu}_i(d\xi ds) \right| ^2\right] \\[2mm] \dsp
= \frac 1{N^2}\sup_{ B\in 2^{\{ w_i; i=1,2,\dots ,N\}}} \max _{k=0,1,\dots ,N-1}E\left[ \sup _{t\in [0,T]}\left| \sum_{i;\, w_i\in B} \int _{s\in (0,t]}\int _{\xi \in [0,\infty )} \right.\right. \\[2mm] \dsp
\left.\left.\phantom{= \frac 1{N^2} \sum_{i;\, w_i\in B} \int _{s\in (0,t]}\int _{\xi \in [0,\infty )}}
\chrfcn{Y_i^{(N)}(s-)\geq Y_C^{(N)}(k/N,0,s-)}\chrfcn{\xi \in [0,w_{i}(Y_i^{(N)}(s-),s))} \tilde{\nu}_i(d\xi ds) \right| ^2\right] .
\arre
\]
Hence, by Doob's martingale inequality (see (6.16) of Chapter I in \cite{IW}) and (3.9) of Chapter II in \cite{IW} we have
\[
\arrb{l} \dsp
E\left[ \sup_{B\in {\mathscr B}(W)} \sup _{y_0 \in [0,1)}\ \sup _{t\in [0,T]}\left| M\ssN_{U}(B,y_0,0,t)\right| ^2\right] \\[2mm] \dsp
\leq \frac {C_3}{N^2} \max _{k=0,1,\dots ,N-1}E\left[ \sum_{i} \int _{s\in (0,T]}\int _{\xi \in [0,\infty )} \chrfcn{\xi\in [0,w_i (Y_i^{(N)}(s-),s))} d\xi ds \right] \\[2mm] \dsp
\leq \frac {C_3R_w(T)T}N 
\arre
\]
where $C_3$ is a positive constant.
Thus, it holds that
\eqnb
\eqna{prfofHDL1007}
\lim _{N\rightarrow \infty} E\left[  \sup_{B\in {\mathscr B}(W)} \sup _{y_0 \in [0,1)}\ \sup _{t\in [0,T]}\left| M\ssN_{U}(B,y_0,0,t)\right| ^2\right] =0 .
\eqne
By \eqnu{MU} again, similarly to the case of $t_0=0$
\[
\arrb{l} \dsp
E\left[ \sup_{B\in {\mathscr B}(W)} \sup _{0\leq t_0 \leq t \leq T}\left| M\ssN_{U}(B,0,t_0,t)\right| ^2\right] \\[2mm] \dsp
= \frac 1{N^2}\sup_{ B\in 2^{\{ w_i; i=1,2,\dots ,N\}}} E\left[ \sup _{0\leq t_0 \leq t \leq T}\left| \sum_{i;\, w_i\in B} \int _{s\in (t_0,t]}\int _{\xi \in [0,\infty )} \right. \right. \\[2mm] \dsp
\left.\left.\phantom{= \frac 1{N^2} \sum_{i;\, w_i\in B} \int _{s\in (0,t]}\int _{\xi \in [0,\infty )}}
\chrfcn{Y_i^{(N)}(s-)\geq Y_C^{(N)}(0,t_0,s-)} \chrfcn{\xi \in [0,w_i (Y_i^{(N)}(s-),s))} \tilde{\nu}_i(d\xi ds)\right| ^2\right] \\[2mm] \dsp
= \frac 1{N^2} \sup_{ B\in 2^{\{ w_i; i=1,2,\dots ,N\}}} E\left[ \sup _{0\leq t_0 \leq t \leq T}\left| \sum_{i;\, w_i\in B} \int _{s\in (0,t]}\int _{\xi \in [0,\infty )} \right. \right. \\[2mm] \dsp
\left.\left.\phantom{= \frac 1{N^2} \sum_{i;\, w_i\in B} \int _{s\in (0,t]}\int _{\xi \in [0,\infty )}}
\chrfcn{Y_i^{(N)}(s-)\geq Y_C^{(N)}(0,t_0,s-)}\chrfcn{\xi \in[0,w_i (Y_i^{(N)}(s-),s))} \tilde{\nu}_i(d\xi ds) \right. \right. \\[2mm] \dsp
\phantom{ \sum_{i;\, w_i\in B} \int _0^t}\left. \left. - \sum_{i;\, w_i\in B} \int _{s\in (0,t_0]}\int _{\xi \in [0,\infty )} \chrfcn{Y_i^{(N)}(s-)\geq Y_C^{(N)}(0,t_0,s-)} \chrfcn{\xi \in [0,w_i (Y_i^{(N)}(s-),s))} \tilde{\nu}_i(d\xi ds)\right| ^2\right] \\[2mm] \dsp
\leq \frac 1{N^2} \sup_{ B\in 2^{\{ w_i; i=1,2,\dots ,N\}}} E\left[ \left( 2\sup _{0\leq t \leq T}\left| \sum_{i;\, w_i\in B} \int _{s\in (0,t]}\int _{\xi \in [0,\infty )}\right. \right. \right. \\[2mm] \dsp
\left.\left.\left.\phantom{= \frac 1{N^2} \sum_{i;\, w_i\in B} \int _{s\in (0,t]}\int _{\xi \in [0,\infty )}}
\chrfcn{Y_i^{(N)}(s-)\geq Y_C^{(N)}(0,t_0,s-)} \chrfcn{\xi \in[0,w_i(Y_i^{(N)}(s-),s))} \tilde{\nu}_i(d\xi ds)\right| \right) ^2\right] .\\[2mm] \dsp
\arre
\]
Hence, by Doob's martingale inequality and (3.9) of Chapter II in \cite{IW} imply
\[
\arrb{l} \dsp
E\left[ \sup_{B\in {\mathscr B}(W)} \sup _{0\leq t_0 \leq t \leq T}\left| M\ssN_{U}(B,0,t_0,t)\right| ^2\right] \\[2mm] \dsp
\leq \frac {4C_3}{N^2} E\left[ \sum_{i} \int _{s\in (0,T]}\int _{\xi \in [0,\infty )} \chrfcn{\xi \in [0,w_i(Y_i^{(N)}(s-),s))} d\xi ds \right] \\[2mm] \dsp
\leq \frac {4C_3R_w(T)T}N.
\arre
\]
Thus, we obtain
\eqnb
\eqna{prfofHDL1010}
\lim _{N\rightarrow \infty}E\left[ \sup_{B\in {\mathscr B}(W)} \sup _{0\leq t_0 \leq t \leq T}\left| M\ssN_{U}(B,0,t_0,t)\right| ^2\right] =0 .
\eqne
Combining \eqnu{prfofHDL1006}, \eqnu{prfofHDL1007} and \eqnu{prfofHDL1010}, we have
\eqnb
\eqna{prfofHDL1011}
\lim _{N\rightarrow \infty} E\left[ \sup_{B\in {\mathscr B}(W)} \sup _{(y_0,t_0) \in \Gamma ;\, t_0\leq T}\ \sup _{t\in [t_0,T]}\left| M\ssN_{U}(B,y_0,t_0,t)\right| ^2\right] =0 .
\eqne

On the other hand, by \eqnu{Uinitconvass} we have
\eqnb
\eqna{prfofHDL1003}
\lim _{N\rightarrow \infty }|| U\ssN(\cdot ,0,t_0)-U(\cdot ,0,t_0)||_{\rm var}= \lim _{N\rightarrow \infty }|| U\ssN(\cdot ,0,0)-U_0(\cdot ,0)||_{\rm var} =0.
\eqne
Thus, \eqnu{prfofHDL1003} and \eqnu{Uinitconvass} implies
\eqnb
\eqna{prfofHDL1005}
\lim _{N\rightarrow \infty} \sup _{(y_0, t_0)\in \Gamma} || U\ssN(\cdot ,y_0,t_0)-U(\cdot ,y_0,t_0)||_{\rm{var}} =0.
\eqne

\eqnu{prfofHDL7W}, \eqnu{prfofHDL1011} and \eqnu{prfofHDL1005} yields
\[
\lim _{N\rightarrow \infty} E\left[ \sup _{t\in [0,T]}W^{(N)}(t)^2\right] = 0.
\]
Hence, there exists a subsequence $\{ N(k)\}$ such that
\[
\lim _{k\rightarrow \infty} \sup _{t\in [0,T]}W^{(N(k))}(t) = 0
\]
almost surely. However, the argument above is also available even if we replace $N$ by any subsequence $N(k)$.
Therefore, we have
\[
\lim _{N\rightarrow \infty} \sup _{t\in [0,T]}W^{(N)}(t) = 0
\]
almost surely. This proves the first assertion of \thmu{HDLposdep}.

We turn to a proof of the second assertion of \thmu{HDLposdep}.
First, we show the uniqueness of the stochastic differential equation \eqnu{taggedlimit}.
Note that
\[
E\left[ \int _{s\in (0,t]} \int _{\xi \in [0,R_w(T)]} \nu _i(d\xi ds)\right] = tR_w(T)
\]
and that for all $i$
\[ \arrb{l} \dsp
\int _{s\in (0,t]} \int _{\xi \in [0,\infty )} Y_i(s-) \chrfcn{\xi \in [0,w_i(Y_i(s-),s))} \nu _i(d\xi ds) \\ \dsp
= \int _{s\in (0,t]} \int _{\xi \in [0,R_w(T)]} Y_i(s-) \chrfcn{\xi \in [0,w_i(Y_i(s-),s))} \nu _i(d\xi ds).
\arre \]
Moreover, there exists a constant $C_T$ such that 
\[
\sup _{s\in [0,T]}|V(W,x,s)-V(W,y,s)|\leq C_T|x-y|, \ x,y\in [0,1).
\]
The proof of Theorem 9.1 in Chapter IV of \cite{IW} is available by taking $U:=[0,R_w(T)]$ and $U_0:=\emptyset$.
Thus, we obtain the uniqueness.

Next, we show that $(Y^{(N)}_1(t), Y^{(N)}_2(t), \dots , Y^{(N)}_L(t))$ converges to $(Y_1(t), Y_2(t),\dots , Y_L(t))$ uniformly in $t\in [0,T]$ almost surely and also converges in the sense of $L^2$.
Let $i\in \{ 1,2,\dots ,L\}$ be fixed.
By \eqnu{prfofHDL4} it is easy to see
\eqnb
\eqna{prfofHDL1020}
\arrb{l} \dsp
Y^{(N)}_i(t) = y^{(N)}_i + M^{(N)}_i(t) + \int _0^t V^{(N)}(W,Y_i^{(N)}(s-),s)ds \\[2mm] \dsp
\phantom{Y^{(N)}_i(t) = }- \int _{s\in (0,t]} \int _{\xi \in [0,\infty )} Y_i^{(N)}(s-) \chrfcn{\xi \in [0,w_i(Y_i^{(N)}(s-),s))} \nu _i(d\xi ds)
\arre
\eqne
where
\eqnb
\eqna{prfofHDL1021}
M^{(N)}_i(t):=\frac 1N \sum _{j=1}^N \int _{s\in (0,t]} \int _{\xi \in [0,\infty )} \chrfcn{Y^{(N)}_i(s-) < Y^{(N)}_j(s-)} \chrfcn{\xi \in [0,w_j (Y_j^{(N)}(s-),s))}\tilde \nu _j(d\xi ds).
\eqne
Hence, \eqnu{taggedlimit} and \eqnu{prfofHDL1020} imply
\eqnb
\eqna{prfofHDL1022}
\arrb{l} \dsp
E\left[ \sup _{s\in [0,t]}|Y^{(N)}_i(s) - Y_i(s)|^2\right] \\[2mm] \dsp
\leq 4|y^{(N)}_i -y_i|^2 + 4E\left[ \sup _{s\in [0,t]}|M^{(N)}_i(s)|^2\right] \\[2mm] \dsp
\phantom{\leq} + 4\int _0^t E\left[ |V^{(N)}(W, Y_i^{(N)}(s-),s)- V(W,Y_i(s-),s)|^2\right] ds \\[2mm] \dsp
\phantom{\leq} + 4E\left[  \sup _{s\in [0,t]}\left| \int _{u\in (0,s]} \int _{\xi \in [0,\infty )} [Y_i^{(N)}(u-) \chrfcn{\xi \in [0,w_i(Y_i^{(N)}(u-),u))}\right. \right. \\[2mm] \dsp
\left. \left. \phantom{\leq E[\sup _{s\in [0,t]}| \int _0^s \int _{\xi \in [0,\infty )}}\qquad -Y_i(u-) \chrfcn{\xi \in [0,w_i(Y_i(u-),u))}] \nu _i(d\xi du)\right| ^2\right] \\[2mm] \dsp
\arre
\eqne
By \eqnu{prfofHDL3} and \eqnu{prfofHDL2}, we have
\[
\arrb{l}\dsp
V\ssN(W,Y\ssN _i(t-),t)-V(W,Y_i(t-),t)
\\[2mm] \dsp {}
=
\int _W w(Y\ssN _i(t-),t)\,
\bigl(U\ssN(dw,Y\ssN _i(t-),t)-U(dw,Y_i(t-),t)\bigr)
\\[2mm] \dsp \phantom{=}
+\int _W \bigl(w(Y\ssN _i(t-),t)-w_i(Y_i(t-),t)\bigr)
\,\,U(dw,Y_i(t-),t)
\\[2mm] \dsp \phantom{=}
+\int_{Y\ssN _i(t-)}^1 \int _W \pderiv{w}{z}(z,t)\,
\bigl(U\ssN(dw,z,t)-U(dw,z,t)\bigr)\,dz
\\[2mm] \dsp \phantom{=}
-\int_{Y_i(t-)}^{Y\ssN _i(t-)} \int _W \pderiv{w}{z}(z,t)\,
U(dw,z,t)\,dz.
\arre
\]
Hence, noting that $0\leq U\ssN \leq 1$ and $0\leq U\leq 1$, there exist positive constants $C_5$ and $C_6$ such that
\eqnb
\eqna{prfofHDL1023}
\arrb{l}\dsp
|V\ssN(W,Y\ssN _i(t-),t)-V(W,Y_i(t-),t)|
\\[2mm] \dsp {}
\leq
C_4 |Y\ssN _i(t-)-Y_i(t-)|+C_5 \sup _{z\in [0,1)}||U\ssN(\cdot ,z,t)-U(\cdot ,z,t)||_{\rm{var}}.
\arre
\eqne
Now we show that
\eqnb
\eqna{prfofHDL1101}
\int _{\xi \in [0,\infty )} |x\chrfcn{\xi \in [0,w_i(x,t))}-y\chrfcn{\xi \in[0,w_i (y,t))}|^2 d\xi \leq C_4|x-y|,\quad x,y\in [0,1)
\eqne
where $C_4$ is a positive constant.
Let $x,y\in [0,1)$ and consider the case that $w_i (y,t) \leq w_i (x,t)$.
Then,
\[
\arrb{l} \dsp
\int _{\xi \in [0,\infty )} |x\chrfcn{\xi \in [0,w_i (x,t))}-y\chrfcn{\xi \in [0,w_i (y,t))}|^2 d\xi \\[2mm] \dsp
= \int _{\xi \in [0,\infty )} |(x-y)\chrfcn{\xi \in [0,w_i (x,t))} + y\chrfcn{\xi \in [w_i (y,t),w_i (x,t))} |^2 d\xi \\[2mm] \dsp
\leq \int _{\xi \in [0,\infty )} \left( 2(x-y)^2\chrfcn{\xi \in [0,w_i (x,t))} + 2y^2\chrfcn{\xi \in [w_i (y,t),w_i (x,t))} \right) d\xi \\[4mm] \dsp
= 2(x-y)^2 w_i (x,t) + y^2(w_i (x,t)-w_i (y,t)).
\arre
\]
Since $w_i$ and the spatial derivative of $w_i$ are bounded, we have
\[
\int _{\xi \in [0,\infty )} |x\chrfcn{\xi \in [0,w_i(x,t))}-y\chrfcn{\xi \in [0,w_i(y,t))}|^2 d\xi
\leq C_4|x-y|,\ x,y\in [0,1),
\]
where $C_4$ is a positive constant.
Therefore, \eqnu{prfofHDL1101} holds.
The case that $w_i(y,t) \geq w_i(x,t)$ is shown similarly.
By \eqnu{compensatedPp}, Doob's martingale inequality and (3.9) of Chapter II in \cite{IW}, there exists a positive constant $C_6$ and we have
\[
\arrb{l}\dsp
E\left[  \sup _{s\in [0,t]}\left| \int _{u\in (0,s]} \int _{\xi \in [0,\infty )} [Y_i^{(N)}(u-) \chrfcn{\xi \in [0,w_i(Y_i^{(N)}(u-),u))} -Y_i(u-) \chrfcn{\xi \in [0,w_i(Y_i(u-),u))}] \nu _i(d\xi du)\right| ^2\right] \\[2mm] \dsp
\leq 2E\left[  \sup _{s\in [0,t]}\left| \int _{u\in (0,s]} \int _{\xi \in [0,\infty )} [Y_i^{(N)}(u-) \chrfcn{\xi \in [0,w_i(Y_i^{(N)}(u-),u))} -Y_i(u-) \chrfcn{\xi \in [0,w_i(Y_i(u-),u))}] \tilde \nu _i(d\xi du)\right| ^2\right] \\[2mm] \dsp
\quad + 2E\left[  \sup _{s\in [0,t]}\left| \int _{u\in (0,s]} \int _{\xi \in [0,\infty )} [Y_i^{(N)}(u-) \chrfcn{\xi \in [0,w_i(Y_i^{(N)}(u-),u))} -Y_i(u-) \chrfcn{\xi \in [0,w_i(Y_i(u-),u))}]d\xi du\right| ^2\right] \\[2mm] \dsp
\leq 2C_6 E\left[  \sup _{s\in [0,t]}\int _{u\in (0,s]} \int _{\xi \in [0,\infty )} \left| Y_i^{(N)}(u-) \chrfcn{\xi \in [0,w_i(Y_i^{(N)}(u-),u))} -Y_i(u-) \chrfcn{\xi \in [0,w_i(Y_i(u-),u))}\right| ^2 d\xi du\right] \\[2mm] \dsp
\quad + 2E\left[  \sup _{s\in [0,t]}\left| \int _0^s[Y_i^{(N)}(u-) w_i(Y_i^{(N)}(u-),u) -Y_i(u-) w_i(Y_i(u-),u)] du\right| ^2\right] .
\arre
\]
By \eqnu{prfofHDL1101} and boundedness of the spatial derivative of $w_i$, there exists a positive constant $C_7$ such that
\[
\arrb{l}\dsp
E\left[  \sup _{s\in [0,t]}\left| \int _{u\in (0,s]} \int _{\xi \in [0,\infty )} [Y_i^{(N)}(u-) \chrfcn{\xi \in [0,w_i(Y_i^{(N)}(u-),u))} -Y_i(u-) \chrfcn{\xi \in [0,w_i(Y_i(u-),u))}] \nu _i(d\xi du)\right| ^2\right] \\[2mm] \dsp
\leq 2C_4 C_6 E\left[ \int _0^t\left| Y_i^{(N)}(u-) -Y_i(u-)\right| du\right] + 2C_7 E\left[  \sup _{s\in [0,t]}\left( \int _0^s\left| Y_i^{(N)}(u-) -Y_i(u-)\right| du\right) ^2\right] \\[2mm] \dsp
\leq 2C_4 C_6 \int _0^t E\left[\left| Y_i^{(N)}(u-) -Y_i(u-)\right| ^2 \right] du+ 2C_7 t E\left[ \int _0^t\left| Y_i^{(N)}(u-) -Y_i(u-)\right| ^2du \right] .
\arre
\]
Thus, we obtain for $t\in [0,T]$
\eqnb
\eqna{prfofHDL1024}
\arrb{l}\dsp
E\left[  \sup _{s\in [0,t]}\left| \int _{u\in (0,s]} \int _{\xi \in [0,\infty )} [Y_i^{(N)}(u-) \chrfcn{\xi \in [0,w_i(Y_i^{(N)}(u-),u))}\right. \right. \\[2mm] \dsp
\left. \left. \phantom{\leq E[\sup _{s\in [0,t]}| \int _0^s \int _{\xi \in [0,\infty )}} -Y_i(u-) \chrfcn{\xi \in [0,w_i(Y_i(u-),u))}] \nu _i(d\xi du)\right| ^2\right] \\[2mm] \dsp
\leq (2C_4 C_6 + 4TC_7) \int _0^t E\left[\left| Y_i^{(N)}(u-) -Y_i(u-)\right| ^2\right] du.
\arre
\eqne
Hence, \eqnu{prfofHDL1022}, \eqnu{prfofHDL1023} and \eqnu{prfofHDL1024} imply that for  $t\in [0,T]$
\[
\arrb{l} \dsp
E\left[ \sup _{s\in [0,t]}|Y^{(N)}_i(s) - Y_i(s)|^2\right] \\[2mm] \dsp
\leq 4|y^{(N)}_i -y_i|^2 + 4E\left[ \sup _{s\in [0,t]}|M^{(N)}_i(s)|^2\right] \\[2mm] \dsp
\phantom{\leq} + 4\int _0^t E\left[ ( C_4 |Y\ssN _i(s-)-Y_i(s-)| + C_5 \sup _{z\in [0,1),s \in[0,T]}||U\ssN(\cdot ,z,s)-U(\cdot ,z,s)||_{\rm{var}})^2\right] ds \\[2mm] \dsp
\phantom{\leq} + 4(2C_4 C_6 + 4TC_7) \int _0^t E\left[\left| Y_i^{(N)}(u-) -Y_i(u-)\right| ^2\right] du \\[2mm] \dsp
\leq 4|y^{(N)}_i -y_i|^2 + 4E\left[ \sup _{s\in [0,t]}|M^{(N)}_i(s)|^2\right] +8C_5^2 \sup _{z\in [0,1),s \in[0,T]}||U\ssN(\cdot ,z,s)-U(\cdot ,z,s)||_{\rm{var}}^2 \\[2mm] \dsp
\phantom{\leq} + [8C_4^2+ 4(2C_4 C_6 + 4TC_7) ]\int _0^t E\left[ \sup _{u\in [0,s]}|Y\ssN _i(u)-Y_i(u)| ^2\right] ds. \\[2mm] \dsp
\arre
\]
By Gronwall's inequality, we obtain
\eqnb
\eqna{prfofHDL1025}
\arrb{l} \dsp
E\left[ \sup _{t\in [0,T]}|Y^{(N)}_i(t) - Y_i(t)| ^2\right] \\[2mm] \dsp
\leq 4e^{C_8T}\left( |y^{(N)}_i -y_i|^2 + E\left[ \sup _{t\in [0,T]}|M^{(N)}_i(t)|^2\right] \right. \\[2mm] \dsp
\left. \phantom{4E\left[ \sup _{s\in [0,t]}|M^{(N)}_i(s)|^2\right]} +2C_5^2 \sup _{z\in [0,1),s \in[0,T]}||U\ssN(\cdot ,z,s)-U(\cdot ,z,s)||_{\rm{var}})^2\right)
\arre
\eqne
where $C_8$ is a positive constant.
Doob's martingale inequality and (3.9) of Chapter II in \cite{IW} and \eqnu{prfofHDL1021} imply there exists a positive constant $C_9$ such that
\[
\arrb{l} \dsp
E\left[ \sup _{s\in [0,t]}|M^{(N)}_i(s)|^2\right] \\[2mm] \dsp
\leq \frac {C_9}{N^2} E\left[ \sum _{j=1}^N \int _{s\in (0,t]} \int _{\xi \in [0,\infty )} \chrfcn{ Y^{(N)}_i(s-) < Y^{(N)}_j(s-)} \chrfcn{\xi \in [0,w_j (Y_j^{(N)}(s-),s))} d\xi ds\right] \\[2mm] \dsp
\leq \frac {C_9R_w(T)t}{N}.
\arre
\]
Hence,
\eqnb
\eqna{prfofHDL1026}
\lim _{N\rightarrow \infty} E\left[ \sup _{s\in [0,t]}|M^{(N)}_i(s)|^2\right] =0.
\eqne
Therefore, by the first assertion of \thmu{HDLposdep}, \eqnu{taggptclconvass}, \eqnu{prfofHDL1025} and \eqnu{prfofHDL1026} we obtain
\[
\lim _{N\rightarrow \infty} E\left[ \sup _{t\in [0,T]}|Y^{(N)}_i(t) - Y_i(t)|^2\right] =0
\]
for $i=1,2,\dots ,L$.
This implies that $(Y^{(N)}_1(t), Y^{(N)}_2(t), \dots , Y^{(N)}_L(t))$ converges to $(Y_1(t), Y_2(t),\dots , Y_L(t))$ uniformly in $t\in [0,T]$ in the sense of $L^2$.

To show the almost sure convergence, see that there exists a subsequence $\{ N(k)\}$ such that $(Y^{(N(k))}_1(t), Y^{(N(k))}_2(t), \dots , Y^{(N(k))}_L(t))$ converges to $(Y_1(t), Y_2(t),\dots , Y_L(t))$ uniformly in $t\in [0,T]$ almost surely.
However, the argument above is also available even if we replace $N$ by any subsequence $N(k)$.
Therefore, we have $(Y^{(N)}_1(t), Y^{(N)}_2(t), \dots , Y^{(N)}_L(t))$ converges to $(Y_1(t), Y_2(t),\dots , Y_L(t))$ uniformly in $t\in [0,T]$ almost surely.

\section{Appendix}

\prpb \label{appendix}
Let $\{ \phi _n\}$ be nondecreasing functions on $[0,1]$ and $\phi$ be a continuous function on $[0,1]$.
Assume that $\phi _n(x)$ converges to $\phi (x)$ for all $x\in [0,1]$.
Then, $\phi _n(x)$ converges to $\phi(x)$ uniformly in $x\in [0,1]$.
\DDD\prpe

\prfb
Let $\varepsilon >0$.
Since $\phi$ is uniformly continuous on $[0,1]$, we can choose a positive integer $N$ such that
\[
|\phi (x) -\phi (y)|< \varepsilon , \quad |x-y| \leq \frac{1}{N}.
\]
By the assumption, there exists a integer $n_0$ such that
\[
\left| \phi _n \left( \frac{k}{N}\right) - \phi \left( \frac{k}{N}\right) \right| < \varepsilon, \quad n\geq n_0\ \mbox{and}\ k=1,2,\dots ,N.
\]
For all $x\in [0,1]$ we can choose $k_x \in \{ 1,2, \dots ,N\}$ such that $0\leq x-k_x /N \leq 1/N$.
Hence, we have for all $x\in [0,1]$ and $n\geq n_0$
\begin{eqnarray*}
&& |\phi _n (x) -\phi (x)| \\
&& \leq \left| \phi _n \left( x\right) - \phi _n \left( \frac{k_x}{N}\right) \right| + \left| \phi _n \left( \frac{k_x}{N}\right) - \phi \left( \frac{k_x}{N}\right) \right| +\left| \phi \left( \frac{k_x}{N}\right) - \phi \left( x\right) \right| \\
&& < \phi _n \left( \frac{k_x+1}{N}\right) - \phi _n\left( \frac{k_x}{N}\right) + 2\varepsilon \\
&& \leq \left| \phi _n \left( \frac{k_x+1}{N}\right) - \phi \left( \frac{k_x+1}{N}\right) \right| + \left| \phi \left( \frac{k_x+1}{N}\right) - \phi \left( \frac{k_x}{N}\right) \right| +\left| \phi \left( \frac{k_x}{N}\right) - \phi _n\left( \frac{k_x}{N}\right) \right| + 2\varepsilon \\
&& \leq 5\varepsilon .
\end{eqnarray*}
This completes the proof.
\QED\prfe

\end{document}